\documentclass[a4paper,twosided, 10pt]{amsart}

\usepackage{amsfonts}
\usepackage{amssymb}
\usepackage{latexsym}
\usepackage{graphics}
\usepackage[2emode]{psfrag}
\usepackage{mathrsfs}
\usepackage{amsthm}
\usepackage{amsmath}
\usepackage[all]{xy}
\usepackage{enumerate}
\usepackage[latin1]{inputenc}
\usepackage{lscape}

\xyoption{2cell} \xyoption{curve} \xyoption{rotate} \xyoption{v2}

\begin{document}
\allowdisplaybreaks

\newcommand{\thlabel}[1]{\label{th:#1}}
\newcommand{\thref}[1]{Theorem~\ref{th:#1}}
\newcommand{\selabel}[1]{\label{se:#1}}
\newcommand{\seref}[1]{Section~\ref{se:#1}}
\newcommand{\lelabel}[1]{\label{le:#1}}
\newcommand{\leref}[1]{Lemma~\ref{le:#1}}
\newcommand{\prlabel}[1]{\label{pr:#1}}
\newcommand{\prref}[1]{Proposition~\ref{pr:#1}}
\newcommand{\colabel}[1]{\label{co:#1}}
\newcommand{\coref}[1]{Corollary~\ref{co:#1}}
\newcommand{\relabel}[1]{\label{re:#1}}
\newcommand{\reref}[1]{Remark~\ref{re:#1}}
\newcommand{\exlabel}[1]{\label{ex:#1}}
\newcommand{\exref}[1]{Example~\ref{ex:#1}}
\newcommand{\delabel}[1]{\label{de:#1}}
\newcommand{\deref}[1]{Definition~\ref{de:#1}}
\newcommand{\eqlabel}[1]{\label{eq:#1}}
\newcommand{\equref}[1]{(\ref{eq:#1})}

\newcommand{\norm}[1]{\| #1 \|}
\def\N{\mathbb N}
\def\Z{\mathbb Z}
\def\Q{\mathbb Q}
\def\mod{\textit{\emph{~mod~}}}
\def\R{\mathcal R}
\def\S{\mathcal S}
\def\*C{{^*\mathcal C}}
\def\C{\mathcal C}
\def\D{\mathcal D}
\def\J{\mathcal J}
\def\M{\mathcal M}
\def\T{\mathcal T}

\newcommand{\Hom}{{\rm Hom}}
\newcommand{\End}{{\rm End}}
\newcommand{\Fun}{{\rm Fun}}
\newcommand{\Mor}{{\rm Mor}\,}
\newcommand{\Aut}{{\rm Aut}\,}
\newcommand{\Hopf}{{\rm Hopf}\,}
\newcommand{\Ann}{{\rm Ann}\,}
\newcommand{\Ker}{{\rm Ker}\,}
\newcommand{\Coker}{{\rm Coker}\,}
\newcommand{\im}{{\rm Im}\,}
\newcommand{\coim}{{\rm Coim}\,}
\newcommand{\Trace}{{\rm Trace}\,}
\newcommand{\Char}{{\rm Char}\,}
\newcommand{\Mod}{{\bf mod}}
\newcommand{\Spec}{{\rm Spec}\,}
\newcommand{\Span}{{\rm Span}\,}
\newcommand{\sgn}{{\rm sgn}\,}
\newcommand{\Id}{{\rm Id}\,}
\newcommand{\Com}{{\rm Com}\,}
\newcommand{\codim}{{\rm codim}}
\newcommand{\Mat}{{\rm Mat}}
\newcommand{\can}{{\rm can}}
\newcommand{\sign}{{\rm sign}}
\newcommand{\kar}{{\rm kar}}
\newcommand{\rad}{{\rm rad}}

\def\Ab{\underline{\underline{\rm Ab}}}
\def\lan{\langle}
\def\ran{\rangle}
\def\ot{\otimes}

\def\id{\textrm{{\small 1}\normalsize\!\!1}}
\def\To{{\multimap\!\to}}
\def\bigperp{{\LARGE\textrm{$\perp$}}}
\newcommand{\QED}{\hspace{\stretch{1}}
\makebox[0mm][r]{$\Box$}\\}

\def\AA{{\mathbb A}}
\def\BB{{\mathbb B}}
\def\CC{{\mathbb C}}
\def\DD{{\mathbb D}}
\def\EE{{\mathbb E}}
\def\FF{{\mathbb F}}
\def\GG{{\mathbb G}}
\def\HH{{\mathbb H}}
\def\II{{\mathbb I}}
\def\JJ{{\mathbb J}}
\def\KK{{\mathbb K}}
\def\LL{{\mathbb L}}
\def\MM{{\mathbb M}}
\def\NN{{\mathbb N}}
\def\OO{{\mathbb O}}
\def\PP{{\mathbb P}}
\def\QQ{{\mathbb Q}}
\def\RR{{\mathbb R}}
\def\SS{{\mathbb S}}
\def\TT{{\mathbb T}}
\def\UU{{\mathbb U}}
\def\VV{{\mathbb V}}
\def\WW{{\mathbb W}}
\def\XX{{\mathbb X}}
\def\YY{{\mathbb Y}}
\def\ZZ{{\mathbb Z}}

\def\aa{{\mathfrak A}}
\def\bb{{\mathfrak B}}
\def\cc{{\mathfrak C}}
\def\dd{{\mathfrak D}}
\def\ee{{\mathfrak E}}
\def\ff{{\mathfrak F}}
\def\gg{{\mathfrak G}}
\def\hh{{\mathfrak H}}
\def\ii{{\mathfrak I}}
\def\jj{{\mathfrak J}}
\def\kk{{\mathfrak K}}
\def\ll{{\mathfrak L}}
\def\mm{{\mathfrak M}}
\def\nn{{\mathfrak N}}
\def\oo{{\mathfrak O}}
\def\pp{{\mathfrak P}}
\def\qq{{\mathfrak Q}}
\def\rr{{\mathfrak R}}
\def\ss{{\mathfrak S}}
\def\tt{{\mathfrak T}}
\def\uu{{\mathfrak U}}
\def\vv{{\mathfrak V}}
\def\ww{{\mathfrak W}}
\def\xx{{\mathfrak X}}
\def\yy{{\mathfrak Y}}
\def\zz{{\mathfrak Z}}

\def\aaa{{\mathfrak a}}
\def\bbb{{\mathfrak b}}
\def\ccc{{\mathfrak c}}
\def\ddd{{\mathfrak d}}
\def\eee{{\mathfrak e}}
\def\fff{{\mathfrak f}}
\def\ggg{{\mathfrak g}}
\def\hhh{{\mathfrak h}}
\def\iii{{\mathfrak i}}
\def\jjj{{\mathfrak j}}
\def\kkk{{\mathfrak k}}
\def\lll{{\mathfrak l}}
\def\mmm{{\mathfrak m}}
\def\nnn{{\mathfrak n}}
\def\ooo{{\mathfrak o}}
\def\ppp{{\mathfrak p}}
\def\qqq{{\mathfrak q}}
\def\rrr{{\mathfrak r}}
\def\sss{{\mathfrak s}}
\def\ttt{{\mathfrak t}}
\def\uuu{{\mathfrak u}}
\def\vvv{{\mathfrak v}}
\def\www{{\mathfrak w}}
\def\xxx{{\mathfrak x}}
\def\yyy{{\mathfrak y}}
\def\zzz{{\mathfrak z}}

\newcommand{\aA}{\mathscr{A}}
\newcommand{\bB}{\mathscr{B}}
\newcommand{\cC}{\mathscr{C}}
\newcommand{\dD}{\mathscr{D}}
\newcommand{\eE}{\mathscr{E}}
\newcommand{\fF}{\mathscr{F}}
\newcommand{\gG}{\mathscr{G}}
\newcommand{\hH}{\mathscr{H}}
\newcommand{\iI}{\mathscr{I}}
\newcommand{\jJ}{\mathscr{J}}
\newcommand{\kK}{\mathscr{K}}
\newcommand{\lL}{\mathscr{L}}
\newcommand{\mM}{\mathscr{M}}
\newcommand{\nN}{\mathscr{N}}
\newcommand{\oO}{\mathscr{O}}
\newcommand{\pP}{\mathscr{P}}
\newcommand{\qQ}{\mathscr{Q}}
\newcommand{\rR}{\mathscr{R}}
\newcommand{\sS}{\mathscr{S}}
\newcommand{\tT}{\mathscr{T}}
\newcommand{\uU}{\mathscr{U}}
\newcommand{\vV}{\mathscr{V}}
\newcommand{\wW}{\mathscr{W}}
\newcommand{\xX}{\mathscr{X}}
\newcommand{\yY}{\mathscr{Y}}
\newcommand{\zZ}{\mathscr{Z}}

\newcommand{\Aa}{\mathcal{A}}
\newcommand{\Bb}{\mathcal{B}}
\newcommand{\Cc}{\mathcal{C}}
\newcommand{\Dd}{\mathcal{D}}
\newcommand{\Ee}{\mathcal{E}}
\newcommand{\Ff}{\mathcal{F}}
\newcommand{\Gg}{\mathcal{G}}
\newcommand{\Hh}{\mathcal{H}}
\newcommand{\Ii}{\mathcal{I}}
\newcommand{\Jj}{\mathcal{J}}
\newcommand{\Kk}{\mathcal{K}}
\newcommand{\Ll}{\mathcal{L}}
\newcommand{\Mm}{\mathcal{M}}
\newcommand{\Nn}{\mathcal{N}}
\newcommand{\Oo}{\mathcal{O}}
\newcommand{\Pp}{\mathcal{P}}
\newcommand{\Qq}{\mathcal{Q}}
\newcommand{\Rr}{\mathcal{R}}
\newcommand{\Ss}{\mathcal{S}}
\newcommand{\Tt}{\mathcal{T}}
\newcommand{\Uu}{\mathcal{U}}
\newcommand{\Vv}{\mathcal{V}}
\newcommand{\Ww}{\mathcal{W}}
\newcommand{\Xx}{\mathcal{X}}
\newcommand{\Yy}{\mathcal{Y}}
\newcommand{\Zz}{\mathcal{Z}}

\numberwithin{equation}{section}
\renewcommand{\theequation}{\thesection.\arabic{equation}}
\newcommand{\bara}[1]{\overline{#1}}
\newcommand{\ContFunt}[2]{\bara{\mathrm{Funt}}(#1,\,#2)}
\newcommand{\lBicomod}[2]{{}_{#1}\mM^{#2}}
\newcommand{\rBicomod}[2]{{}^{#1}\mM_{#2}}
\newcommand{\Bicomod}[2]{{}^{#1}\mM^{#2}}
\newcommand{\lrBicomod}[2]{{}^{#1}\mM^{#2}}
\newcommand{\rcomod}[1]{\mM^{#1}}
\newcommand{\rmod}[1]{\mM_{#1}}
\newcommand{\lmod}[1]{{}_{#1}\mM}
\newcommand{\Bimod}[2]{{}_{#1}\mM_{#2}}
\newcommand{\Sf}[1]{\mathsf{#1}}
\newcommand{\Sof}[1]{S^{{\bf #1}}}
\newcommand{\Scr}[1]{\mathscr{#1}}
\newcommand{\Tof}[1]{T^{{\bf #1}}}
\renewcommand{\hom}[3]{\mathrm{Hom}_{#1}\left(#2,\,#3\right)}
\newcommand{\Coint}[2]{\mathrm{Coint}(#1,#2)}
\newcommand{\IntCoint}[2]{\mathrm{InCoint}(#1,#2)}
\newcommand{\Coder}[2]{\mathrm{Coder}(#1,#2)}
\newcommand{\IntCoder}[2]{\mathrm{InCoder}(#1,#2)}
\newcommand{\lr}[1]{\left(\underset{}{} #1 \right)}
\newcommand{\Ext}[3]{\mathrm{Ext}_{\eE}^{#1}\lr{#2,\,#3}}
\newcommand{\equalizerk}[2]{\mathfrak{eq}_{#1,\,#2}^k}
\newcommand{\equalizer}[2]{\mathfrak{eq}_{#1,\,#2}}
\newcommand{\coring}[1]{\mathfrak{#1}}
\newcommand{\tensor}[1]{\otimes_{#1}}

\def\units{{\mathbb G}_m}
\def\rightact{\hbox{$\leftharpoonup$}}
\def\leftact{\hbox{$\rightharpoonup$}}

\def\*C{{}^*\hspace*{-1pt}{\Cc}}

\def\text#1{{\rm {\rm #1}}}

\def\smashco{\mathrel>\joinrel\mathrel\triangleleft}
\def\cosmash{\mathrel\triangleright\joinrel\mathrel<}

\def\ol{\overline}
\def\ul{\underline}
\def\dul#1{\underline{\underline{#1}}}
\def\Nat{\dul{\rm Nat}}
\def\Set{\dul{\rm Set}}

\renewcommand{\subjclassname}{\textup{2000} Mathematics Subject Classification}

\newtheorem{prop}{Proposition}[section]
\newtheorem{lemma}[prop]{Lemma}
\newtheorem{cor}[prop]{Corollary}
\newtheorem{theo}[prop]{Theorem}

\theoremstyle{definition}
\newtheorem{Def}[prop]{Definition}
\newtheorem{ex}[prop]{Example}
\newtheorem{exs}[prop]{Examples}

\theoremstyle{remark}
\newtheorem{rems}[prop]{Remarks}
\newtheorem{rem}[prop]{Remark}

\title[Wide Morita Contexts in Bicategories.]
{Wide Morita Contexts in Bicategories.}
\date{\today}
\author{L. El Kaoutit}
\address{Departamento de \'Algebra. Facultad de Educaci\'on y Humanidades de Ceuta.
Universidad de Granada. El Greco N. 10. E-51002 Ceuta, Spain}
\email{kaoutit@ugr.es}

\keywords{ Bimodules. Corings. Category of comodules. Morita
contexts. Bicategories.\\
Research supported by  grant MTM2007-61673 from the Ministerio de
Educaci\'{o}n y Ciencia of Spain, and P06-FQM-01889 from Junta de Andaluc\'ia.} 

\baselineskip 10pt

\begin{abstract}
We give a formal concept of (right) wide Morita context between two
$0$-cells in an arbitrary bicategory. We then construct a new
bicategory with the same $0$-cells as the previous one, and with
$1$-cells all these (right) wide Morita contexts. An application to
the (right) Eilenberg-Moore bicategory of comonads associated to the
bicategory of bimodules is also given.
\end{abstract}

\maketitle

\section*{Introduction}

A Morita context (see \cite{Bass:1968}) connecting two associative rings with unit $A$ and $B$, is a four-tuple $({}_BM_A,{}_AN_B,\varphi,\psi)$ consisting of two bimodules ${}_BM_A$ and ${}_AN_B$ and two bilinear maps $\varphi: M\tensor{A}N \to B$ and $\psi: N\tensor{B}M \to A$ satisfying compatibility conditions. A morphism between two Morita contexts connecting $A$ and $B$ is a pair of bilinear maps satisfying two equations. A basic result is that if $\varphi$ and $\psi$ are surjective, then they are in fact invertible, and so $A$ and $B$ have equivalent categories of modules. Part of the fundamental Morita theorem says that there is a one-to-one correspondence between the isomorphism types of category equivalences between categories of modules; and the isomorphism types of Morita contexts with surjective maps. A natural question that was posed was how far is a Morita context with not necessarily invertible maps from an equivalence of categories. An answer was given by B. J. M\"{u}ller in \cite{Muller:1974} (based on the works of T. Kato  \cite{Kato:1970, Kato:1973}), which says that every Morita context induces an equivalence between the quotient categories of modules categories determined by the two trace ideals. Part of M\"{u}ller's Theorem says that there is a one-to-one correspondence between the isomorphism types of maximal category equivalences between full
subcategories of categories of modules containing the modules $A_A$ and $B_B$, and the isomorphism types of right normalized Morita contexts connecting $A$ and $B$ (see \cite[Theorem 7]{Muller:1974}). A companion to M\"{u}ller's Theorem was given by W. K. Nicholson and J. F. Watters in \cite[Theorem 5]{Nicholson/Watters:1988}.

The notion of a Morita context was extended to an arbitrary Grothendieck category by
F. Casta\~{n}o Iglesias and J. G\'omez-Torrecillas in \cite[p. 602]{Castano/Gomez:1995}. Let $\Aa$ and $\Bb$ be additive categories with cokernels and direct sums, consider a pair of additive covariant functors $\xymatrix{\Ff:\Bb \ar@<0,5ex>[r] & \Aa:\Gg \ar@<0,5ex>[l]}$
together with two natural transformations $\eta:\Ff\Gg \to \id_{\Aa}$ and $\rho:\Gg\Ff \to \id_{\Bb}$ such that $\Ff,\Gg$ preserve cokernels and $\eta_{\Ff(Y)} \,=\, \Ff \rho_Y$, $\rho_{\Gg(X)}\,=\, \Gg\eta_X$, for every pair of objects $Y \in \Bb$ and $X \in \Aa$. The four-tuple $(\Ff,\Gg,\eta,\rho)$ is referred to as a
\emph{wide right Morita context} between $\Aa$ and $\Bb$. A \emph{wide left Morita context} is dually defined, that is, in the dual categories $\Aa^0$ and $\Bb^0$. Obviously if one of these natural transformations is invertible, then  the context is nothing but an adjunction (with invertible unit or counit) between the categories $\Aa$ and $\Bb$. This is the case, which corresponds to a wide left Morita context, when for example $\Bb$ is a Grothendieck category and $\Ff$ is a localizing functor.

In a series of papers \cite{Castano/Gomez:1995, Castano/Gomez:1996,
Castano/Gomez:1998}, F. Casta\~{n}o Iglesias and J. G\'omez-Torrecillas
proved several results concerning an equivalences of categories
between full subcategories (quotient categories) that can be
constructed form wide (right) Morita contexts. Stimulating
applications to the categories of modules and graded modules were
explicitly exhibited in the first two papers. In the third one
Morita-Takeuchi contexts \cite{Takeuchi:1977} as wide left Morita
contexts between categories of comodules over coalgebras over a
field were particularly studied. A different approach of this case
was also given by I. Berbec in \cite{Berbec:2003}. Recently in
\cite{Chifan/Dascalescu/Nastasescu:2005} N. Chifan, S. D$\rm \breve{a}sc\breve{a}lescu$ and C. N$\rm
\breve{a}st\breve{a}sescu$ unify all
the previous results in the framework of abelian categories where
they also added new examples, see \cite[Section
5]{Chifan/Dascalescu/Nastasescu:2005}. Morita contexts with
injective bilinear maps were studied more recently by J. Y. Abuhlail
and S. K. Nauman in \cite{Abuhlail/Nauman:2008}, where several
equivalences of categories were established \cite[Section
5]{Abuhlail/Nauman:2008}.

Bimodules are $1$-cells in the bicategory $\Sf{Bim}$ (rings,
bimodules, bilinear maps), while functors are $1$-cells in the
$2$-category $\Sf{Cat}$ (categories, functors, natural
transformations). Thus given an arbitrary bicategory (see the
definition below), one can directly extend the concept of wide
right Morita context between any pair of $0$-cells. This means that
the best and unifying setting of studying (wide) Morita contexts
connecting rings or additive categories is the framework of
bicategories. On the other hand, we sincerely think that this
concept can be studied in arbitrary bicategory by its own right as
was done for the notions of adjunction and of equivalence in
bicategories.

In this paper we consider wide right Morita contexts in an arbitrary
bicategory $\Sf{B}$. In the first Section we study the relationship
between internal equivalences and wide right Morita contexts
(Propositions \ref{Int-Eq} and \ref{abelian}). We propose a notion of morphism
of right Morita contexts, and construct a new bicategory which we
denote by $\Sf{W(B)}$. The $0$-cells of $\Sf{W(B)}$ are all that of
$\Sf{B}$ and its Hom-Categories are wide right Morita contexts and
their morphisms (Proposition \ref{Bicategory}). We show that this construction is in fact functorial with respect to homomorphisms of bicategories (Proposition \ref{Morph-bicatg}).  In the second
Section we give an application of this construction taking
$\Sf{B}\,=\,\Sf{REM(Bim)}$ the right Eilenberg-Moore bicategory of
comonad attached to $\Sf{Bim}$, whose $0$-cells are all corings over
rings with identity. This bicategory was studied in
\cite{Brzezinski/Kaoutit/Gomez:2006} and earlier introduced for
general setting in \cite{Lack/Street:2002}. We prove that every
$1$-cell in $\Sf{W(REM(Bim))}$ (i.e. a wide right Morita context)
induces a wide right Morita context between the categories of right
comodules over corings in the sense of \cite{Castano/Gomez:1995}
(Proposition \ref{Morita-comod}). Notice that the categories of
comodules over corings are not in general abelian
\cite{Kaoutit/Gomez/Lobillo:2004c}. The converse of that proposition
is also analyzed (Proposition \ref{Reciproco}).

\textsc{Notations and basic notions.} Since we will construct
bicategories from other one, it is convenient to include the
general definition of a bicategory, more details can be found in
the fundamental paper \cite{Benabou:1967} (see also
\cite{Leinster:1998} for a basic definitions). A \emph{bicategory}
$\Sf{B}$ consists of the following data subject to the
forthcoming
axioms\\
\textsl{Data}:
\begin{enumerate}[$\bullet$]
\item \emph{Objects (or $0$-cells)}. A class of objects $\textsl{Ob}(\Sf{B})$ denoted by
$A,B,C,D, \cdots$.

\item \emph{Hom-Categories}. For each pair of $0$-cells $A$ and
$B$, there is a category ${}_A\Sf{B}_B$. The objects of this
category are known as $1$-cells from $B$ to $A$, and morphisms as
$2$-cells from $B$ to $A$.

\item \emph{Functors}. Defining the horizontal and vertical two
sided multiplications $$\xymatrix@R=0pt{\ccc_{ABC}: {}_A\Sf{B}_B
\times
{}_B\Sf{B}_C \ar@{->}[r] & {}_A\Sf{B}_C \\ (f,g) \ar@{->}[r] & fg \\
(\alpha,\beta) \ar@{->}[r] & \alpha \beta }$$ and $\mathbb{I}:
{\bf 1} \to {}_A\Sf{B}_A$ defining an identity $1$-cell
$\mathbb{I}_A$, where $\bf{1}$ denotes the category with one
object.

\item \emph{Natural isomorphism}. Defining associativity up to
isomorphisms, and compatibility with left and right
multiplications by identity $1$-cells
$$
\xy *+{{}_A\Sf{B}_C \times {}_C\Sf{B}_D}="p",
p+<0pt,2cm>*+{{}_A\Sf{B}_B \times {}_B\Sf{B}_C \times
{}_C\Sf{B}_D}="2", p+<5cm,0pt>*+{{}_A\Sf{B}_D}="3",
p+<5cm,2cm>*+{{}_A\Sf{B}_B \times {}_B\Sf{B}_D}="4",
p+<1.4cm,.8cm>*+{}="5", p+<2cm,1.1cm>*+{}="6", {"2" \ar@{->}_{\ccc
\times 1} "p"}, {"2" \ar@{->}^-{1 \times \ccc} "4"}, {"4"
\ar@{->}^-{\ccc} "3"}, {"p" \ar@{->}_-{\ccc} "3"}, {"5"
\ar@{=>}_-{\aaa} "6"}
\endxy
$$
$$
 \xy *+{{}_A\Sf{B}_B \times {}_B\Sf{B}_B}="p",
p+<0pt,2cm>*+{{}_A\Sf{B}_B \times {\bf 1}}="1",
p+<4cm,0pt>*+{{}_A\Sf{B}_B}="2", p+<.5cm,.5cm>*+{}="3",
p+<1cm,1cm>*+{}="4", {"p" \ar@{->}^-{1 \times \II} "1"}, {"p"
\ar@{->}_-{\ccc} "2"}, {"1" \ar@{->}^-{} "2"}, {"3"
\ar@{=>}_-{\rrr} "4"}
\endxy, \qquad \xy *+{{}_A\Sf{B}_B}="p", p+<4cm,0pt>*+{{}_A\Sf{B}_A
\times {}_A\Sf{B}_B}="1", p+<4cm,2cm>*+{ {\bf 1} \times
{}_A\Sf{B}_B}="2", p+<3cm,.6cm>*{}="3",p+<2.5cm,1cm>*{}="4", {"3"
\ar@{=>}^-{\lll} "4"}, {"2" \ar@{->}_-{} "p"}, {"2" \ar@{->}^-{\II
\times 1} "1"}, {"1" \ar@{->}^-{\ccc} "p"}
\endxy
$$
that is, $2$-cells $$\xymatrix{ \aaa_{fgh}: \,\,
(fg)h \ar@{->}^-{\cong}[r] & f(gh), & \rrr_f:\,\, f \II_B
\ar@{->}^-{\cong}[r] & f, & \lll_f:\,\, \II_A f
\ar@{->}^-{\cong}[r] & f.}$$
\end{enumerate}
\textsl{Axioms}:\\ The following diagrams commute
$$
\xy *+{((fg)h)e}="p", p+<2.5cm,1.5cm>*+{(f(gh))e}="1",
p+<5cm,1.5cm>*+{f((gh)e)}="2", p+<7.5cm,0pt>*+{f(g(he))}="3",
p+<3.5cm,-1.5cm>*+{(fg)(he)}="4", p+<3.3cm,0pt>*+{\text{Ax}\; \rm
I}, {"p" \ar@{->}^-{\aaa 1} "1"}, {"p" \ar@{->}_-{\aaa} "4"}, {"1"
\ar@{->}^-{\aaa} "2"}, {"2" \ar@{->}^-{1 \aaa} "3"}, {"4"
\ar@{->}_-{\aaa} "3"}, p+<10cm,1.5cm>*+{(f\II_B)g}="7",
p+<13cm,1.5cm>*+{f(\II_B g)}="8", p+<11.5cm, -1.5cm>*+{fg}="9",
p+<11.5cm, .3cm>*+{\text{Ax}\; \rm II}, {"7" \ar@{->}^-{\aaa}
"8"}, {"7" \ar@{->}_-{\rrr 1} "9"}, {"8" \ar@{->}^-{1\lll} "9"}
\endxy
$$

Using Axioms I and II, and the naturality of $\aaa$, $\rrr$ and
$\lll$ we obtain the equality $\rrr_{\II_{-}}=\lll_{\II_{-}}$, and
that the following diagrams are commutative
\begin{equation}\label{AxIII}
\xy *+{(fg)\II_C}="p", p+<3cm, 0pt>*+{f(g\II_C)}="1",
p+<1.5cm,-2.5cm>*+{fg}="2", {"p" \ar@{->}^-{\aaa} "1"},
{"p"\ar@{->}_-{\rrr} "2"}, {"1" \ar@{->}^-{1\rrr} "2"}
\endxy\qquad \xy *+{\II_A(fg)}="p", p+<3cm, 0pt>*+{(\II_Af)g}="1",
p+<1.5cm,-2.5cm>*+{fg}="2", {"p" \ar@{->}^-{\aaa^{-1}} "1"},
{"p"\ar@{->}_-{\lll} "2"}, {"1" \ar@{->}^-{\lll 1} "2"}
\endxy
\end{equation}
\begin{equation}\label{AxIV}
\xy *+{(\II_Af)\II_B}="p", p+<3cm, 0pt>*+{\II_A(f\II_B)}="1",
p+<0pt,-2.5cm>*+{f\II_B}="2", p+<3cm, -2.5cm>*+{\II_Af}="3",
p+<1.5cm, -3.5cm>*+{f}="4", {"p" \ar@{->}^-{\aaa} "1"},
{"p"\ar@{->}_-{\lll 1} "2"}, {"1" \ar@{->}^-{1\rrr} "3"}, {"2"
\ar@{->}_-{\rrr} "4"}, {"3" \ar@{->}^-{\lll} "4"}
\endxy
\end{equation}

A \emph{morphisms} $(\Scr{F}, \Phi)$ from a bicategory $(\Sf{B},\ccc, \aaa,\rrr,\lll,\II)$ to $(\Sf{B}',\ccc',\aaa',\rrr',\lll',\II')$ consists of the following data subject to the forthcoming axioms
\\
\textsl{Data}:
\begin{enumerate}[$\bullet$]
\item \emph{Function}. $\Scr{F}: \textsl{Ob}(\Sf{B}) \to \textsl{Ob}(\Sf{B}')$.

\item \emph{Functors}. $\Scr{F}_{A,\,B}: {}_A\Sf{B}_B \to {}_{\Scr{F}(A)}\Sf{B'}_{\Scr{F}(B)}$.

\item \emph{Natural isomorphism}.
$$
\xy *+{{}_{\Scr{F}(A)}\Sf{B'}_{\Scr{F}(B)} \times {}_{\Scr{F}(B)}\Sf{B'}_{\Scr{F}(C)}}="p",
p+<0pt,2cm>*+{{}_A\Sf{B}_B \times {}_B\Sf{B}_C}="2", p+<6cm,0pt>*+{{}_{\Scr{F}(A)}\Sf{B'}_{\Scr{F}(C)}}="3",
p+<6cm,2cm>*+{{}_A\Sf{B}_C}="4",
p+<1.3cm,.8cm>*+{}="5", p+<2cm,1.1cm>*+{}="6", {"2" \ar@{->}_{(\Scr{F}_{A,B},\,\Scr{F}_{B,C})} "p"}, {"2" \ar@{->}^-{\ccc_{ABC}} "4"}, {"4"
\ar@{->}^-{\Scr{F}_{A,C}} "3"}, {"p" \ar@{->}_-{\ccc'_{\Scr{F}(A)\Scr{F}(B)\Scr{F}(C)}} "3"}, {"5"
\ar@{=>}_-{\Phi^1} "6"}
\endxy \qquad
\xy *+{{\bf 1}}="p",
p+<0pt,2cm>*+{{\bf 1}}="2", p+<4cm,0pt>*+{{}_{\Scr{F}(A)}\Sf{B'}_{\Scr{F}(A)}}="3",
p+<4cm,2cm>*+{{}_A\Sf{B}_A}="4",
p+<1.3cm,.8cm>*+{}="5", p+<2cm,1.1cm>*+{}="6", {"2" \ar@{=}_{} "p"}, {"2" \ar@{->}^-{\II_A} "4"}, {"4"
\ar@{->}^-{\Scr{F}_{A,A}} "3"}, {"p" \ar@{->}_-{\II'_{\Scr{F}(A)}} "3"}, {"5"
\ar@{=>}_-{\Phi^0} "6"}
\endxy
$$
that is, two $2$-cells
$$\Phi^1_{f,g}: \Scr{F}(f) \Scr{F}(g) \longrightarrow \Scr{F}(fg) \quad \text{ and } \quad \Phi^0_A: \II'_{\Scr{F}(A)} \longrightarrow \Scr{F}(\II_A).
$$
\end{enumerate}
\textsl{Axioms}:\\ The following diagrams commute
$$
\xy *+{\Scr{F}(f)\Scr{F}(gh)}="p", p+<6cm,0pt>*+{\Scr{F}((fg)h)}="1",
p+<6cm,2cm>*+{\Scr{F}(fg)\Scr{F}(h)}="2", p+<0pt, 2cm>*+{\Scr{F}(f)(\Scr{F}(g)\Scr{F}(h))}="3",
p+<2cm,-1.5cm>*+{\Scr{F}(f(gh))}="4", p+<2cm,3.5cm>*+{(\Scr{F}(f)\Scr{F}(g))\Scr{F}(h)}="5",  {"3" \ar@{->}_-{1\Phi} "p"}, {"p" \ar@{->}_-{\Phi} "4"}, {"1"\ar@{->}^-{\Scr{F}(\aaa)} "4"}, {"2" \ar@{->}^-{\Phi} "1"}, {"5"
\ar@{->}^-{\Phi 1} "2"}, {"5" \ar@{->}_-{\aaa'} "3"}
\endxy
$$
$$
\xy *+{\Scr{F}(f)}="p",
p+<0pt,2cm>*+{\Scr{F}(f)\II'_{\Scr{F}(B)}}="2", p+<4cm,0pt>*+{\Scr{F}(f\II_B)}="3",
p+<4cm,2cm>*+{\Scr{F}(f)\Scr{F}(\II_B)}="4",
p+<1.3cm,.8cm>*+{}="5", {"2" \ar@{->}_{\rrr'} "p"}, {"2" \ar@{->}^-{1\Phi} "4"}, {"4"
\ar@{->}^-{\Phi} "3"}, {"p" \ar@{->}_-{{\Scr{F}(\rrr^{-1})}} "3"}
\endxy \qquad \qquad
\xy *+{\Scr{F}(f)}="p",
p+<0pt,2cm>*+{\II'_{\Scr{F}(A)}\Scr{F}(f)}="2", p+<4cm,0pt>*+{\Scr{F}(\II_A f)}="3",
p+<4cm,2cm>*+{\Scr{F}(\II_A)\Scr{F}(f)}="4",
p+<1.3cm,.8cm>*+{}="5", {"2" \ar@{->}_{\lll'} "p"}, {"2" \ar@{->}^-{\Phi 1} "4"}, {"4"
\ar@{->}^-{\Phi} "3"}, {"p" \ar@{->}_-{{\Scr{F}(\lll^{-1})}} "3"}
\endxy
$$
$(\Scr{F}, \Phi)$ is called a \emph{homomorphism} provided that $\Phi$ are natural isomorphisms, that is, $\Scr{F}(f)\Scr{F}(g) \cong \Scr{F}(fg)$ and $\Scr{F}(\II_A) \cong \II'_{\Scr{F}(A)}$.

\section{Wide Morita contexts as $1$-cells}\label{Sect-1}

Let $\Sf{B}$ be a bicategory with natural isomorphisms $\aaa$,
$\rrr$ and $\lll$. Given two $0$-cells $A$ and $B$, we define a
\emph{wide right Morita context} (WRMC for short) from $B$ to $A$
as a four-tuple $\Gamma=(f,g,\eta,\rho)$ consisting of two
$1$-cells $f \in {}_A\Sf{B}_B$ and $g \in {}_B\Sf{B}_A$, and two
$2$-cells $\eta: fg \to \II_A$ and $\rho: gf \to \II_B$ such that
the following diagrams commute
\begin{equation}\label{Wide}
\xy *+{g(fg)}="p", p+<2.5cm,0pt>*+{g\II_A}="1",
p+<4cm,-1cm>*+{g}="2", p+<0pt,-2cm>*+{(gf)g}="3",
p+<2.5cm,-2cm>*+{\II_B g}="4", {"p" \ar@{->}^-{1 \eta} "1"}, {"1"
\ar@{->}^-{\rrr } "2"}, {"3" \ar@{->}^-{\aaa} "p"}, {"3"
\ar@{->}_-{\rho 1} "4"}, {"4" \ar@{->}_-{\lll} "2"}
\endxy \qquad \qquad \xy *+{(fg)f}="p", p+<2.5cm,0pt>*+{\II_A
f}="1", p+<4cm,-1cm>*+{f}="2", p+<0pt,-2cm>*+{f(gf)}="3",
p+<2.5cm,-2cm>*+{f\II_B }="4" , {"p" \ar@{->}^-{\eta 1} "1"}, {"1"
\ar@{->}^-{\lll} "2"}, {"p" \ar@{->}_-{\aaa} "3"}, {"3"
\ar@{->}_-{1 \rho} "4"}, {"4" \ar@{->}_-{\rrr} "2"}
\endxy
\end{equation}
Wide left Morita contexts are dually defined in the conjugate
category $\Sf{B}^{co}$ obtained from $\Sf{B}$ by reversing
$2$-cells.

Let $A$ and $B$ two $0$-cells of $\Sf{B}$. Recall (see
\cite{Leinster:1998}) that $A$ is said to be \emph{equivalent} to
$B$ (\emph{inside} $\Sf{B}$), if there exists a pair of $1$-cells $f
\in {}_A\Sf{B}_B$ and $g \in {}_B\Sf{B}_A$ together with two
isomorphisms $\eta: fg \to \II_A$ in ${}_A\Sf{B}_A$, and
$\theta:\II_B \to gf$ in ${}_B\Sf{B}_B$.

It is well known that any equivalence inside a bicategory induces an adjunction with invertible unit (or counit). Whence it induces a wide (right) Morita context. For sake of completeness, we give a complet proof of this fact.

\begin{prop}\label{Int-Eq}
Assume that $A$ is equivalent to $B$ inside $\Sf{B}$ via the
isomorphisms $\eta: fg \to \II_A$ and $\theta:\II_B \to gf$.
\begin{enumerate}[(i)]
\item The following diagrams commute
$$
\xy *+{(fg)(fg)}="p", p+<3.5cm, 0pt>*+{(fg)\II_A}="1",
p+<0pt,-2cm>*+{\II_A(fg)}="2", p+<3.5cm,-2cm>*+{fg}="4", {"p"
\ar@{->}^-{1\eta} "1"}, {"1" \ar@{->}^-{\rrr} "4"}, {"p"
\ar@{->}_-{\eta 1} "2"}, {"2" \ar@{->}^-{\lll} "4"}
\endxy \qquad \xy *+{\II_B(gf)}="p", p+<3.5cm, 0pt>*+{(gf)(gf)}="1",
p+<0pt,-2cm>*+{gf}="2", p+<3.5cm,-2cm>*+{(gf)\II_B}="4", {"p"
\ar@{->}^-{\theta 1} "1"}, {"2" \ar@{->}^-{\lll^{-1}} "p"}, {"2"
\ar@{->}^-{\rrr^{-1}} "4"}, {"4" \ar@{->}_-{1 \theta} "1"}
\endxy
$$

\item The functor $f-:{}_B\Sf{B}_B \to {}_A\Sf{B}_B$ (respectively
$g-:{}_A\Sf{B}_A \to {}_B\Sf{B}_A$), defined by left
multiplication by $f$ (respectively by $g$), is fully faithful.

\item There exists a $2$-cell $\rho: gf \to \II_B$ such that
$(f,g,\eta,\rho)$ is WRMC from $B$ to $A$.
\end{enumerate}
\end{prop}
\begin{proof}
$(i)$. This is immediate since $\rrr$, $\lll$ are natural and
$\eta$, $\theta$ are monomorphisms.

$(ii)$. Let $\sigma, \sigma': g_1 \to g_2$ two morphisms in the
category  ${}_B\Sf{B}_A$. We claim that, if
$f\sigma\,=\,f\sigma'$, then $\sigma=\sigma'$, that is the left
multiplication by $f$ is faithful. This in fact follows form the
following commutative diagram
$$\xymatrix@R=60pt@C=60pt{  g_1 \ar@{->}^-{\lll^{-1}}[r] \ar@<-0,5ex>_-{\sigma}[d] \ar@<0,5ex>^-{\sigma'}[d]
 & \II_Bg_1 \ar@{->}^-{\theta 1}[r]
\ar@<-0,5ex>_-{ 1\sigma}[d] \ar@<0,5ex>^-{ 1\sigma'}[d] & (gf)g_1
\ar@{->}^-{\aaa}[r] \ar@<-0,5ex>_-{1 \sigma}[d] \ar@<0,5ex>^-{1
\sigma'}[d] & g(fg_1) \ar@<-0,5ex>_-{1(1\sigma)}[d]
\ar@<0,5ex>^-{1(1\sigma')}[d]
\\ g_2 \ar@{->}_-{\lll^{-1}}[r] & \II_Bg_2 \ar@{->}_-{\theta 1}[r]
& (gf)g_2 \ar@{->}_-{\aaa}[r] & g(fg_2)}
$$
whose horizontal maps are all isomorphisms. Similar arguments are
used for the left multiplication by $g$. Consider now a morphism
$\alpha:fh_1 \to fh_2$ in the category ${}_A\Sf{B}_B$ where $h_i$
are objects of the category ${}_B\Sf{B}_B$, for $i=1,2$. So, we
can define the composition
$$
\xy *+{\sigma:\,\, h_1}="p", p+<2cm,0pt>*+{\II_Bh_1}="1",
p+<4cm,0pt>*+{(gf)h_1}="2", p+<6cm,0pt>*+{g(fh_1)}="3",
p+<8cm,0pt>*+{g(fh_2)}="4", p+<10cm,0pt>*+{(gf)h_2}="5",
p+<13cm,0pt>*+{\II_B h_2}="6", p+<14.5cm,0pt>*+{h_2}="7", {"p"
\ar@{->}^-{\lll^{-1}} "1"}, {"1" \ar@{->}^-{\theta 1} "2"}, {"2"
\ar@{->}^-{\aaa} "3"}, {"3" \ar@{->}^-{1\alpha} "4"}, {"4"
\ar@{->}^-{\aaa^{-1}} "5"}, {"5" \ar@{->}^-{\theta^{-1}1} "6"},
{"6" \ar@{->}^-{\lll} "7"}
\endxy
$$
Using the diagrams stated in item $(i)$ and in \eqref{AxIII},
Axioms II and II, in conjunction with the naturality of $\aaa$, we
can check that $g(f\sigma)\,=\,g\alpha$, and so
$f\sigma\,=\,\alpha$ since $g-$ is faithful. This proves that the
left multiplication by $f$ is full. Similar arguments are used for
left multiplication by $g$.

$(iii)$. Consider the following morphism in the category
${}_A\Sf{B}_B$
$$
\xy *+{\kappa:\,\,f(gf)}="p",
p+<3cm,0pt>*+{(fg)f}="1", p+<6cm,0pt>*+{\II_Af}="2",
p+<9cm,0pt>*+{f}="3", p+<12cm,0pt>*+{f\II_B}="4", {"p"
\ar@{->}^-{\aaa^{-1}} "1"}, {"1" \ar@{->}^-{\eta 1} "2"}, {"2"
\ar@{->}^-{\lll} "3"}, {"3" \ar@{->}^-{\rrr^{-1}} "4"}
\endxy
$$
by item $(ii)$ there exists $\rho: gf \to \II_B$ such that
$$f\rho\,\,=\,\,\kappa\,\,=\,\, \rrr^{-1}_f \circ \lll_f \circ \eta f \circ
\aaa_{fgf}^{-1},$$ that is $\rrr_f \circ f\rho \circ
\aaa_{fgf}\,=\, \lll_f \circ \eta f $. This last equality is exactly
the commutativity of the first diagram in \eqref{Wide}. To show
the commutativity of the second diagram in \eqref{Wide}, we
multiply on the left by $f$ to get
$$
\xy *+{f((gf)g)}="p", p+<8cm,0pt>*+{f(g(fg))}="1",
p+<8cm,-4cm>*+{f(g\II_A)}="3", p+<0pt, -4cm>*+{f(\II_Bg)}="2",
p+<4cm,-6cm>*+{(fg)}="4", p+<4cm,-1cm>*+{(fg)(fg)}="2'",
p+<4cm,-2.5cm>*+{((fg)f)g}="3'", p+<4cm,-4cm>*+{(\II_Af)g}="4'",
{"p" \ar@{->}_-{1(\rho 1)} "2"}, {"p" \ar@{->}^-{1\aaa} "1"}, {"1"
\ar@{->}^-{1(1\eta)} "3"}, {"2" \ar@{->}^-{1\lll} "4"}, {"3"
\ar@{->}_-{1\rrr} "4"}, {"1" \ar@{->}_-{\aaa^{-1}} "2'"}, {"2'"
\ar@{->}_-{\aaa^{-1}} "3'"}, {"3'" \ar@{->}_-{(\eta 1)1} "4'"},
{"4'" \ar@{->}_-{\lll 1} "4"}, p+<2cm, -3cm>*+{1}*\cir{}, p+<6cm,
-3cm>*+{2}*\cir{}
\endxy
$$
Diagram $\xy *+{1}*\cir{} \endxy$ commutes by applying consecutively
Axiom I, the naturality of $\aaa$, the equality satisfied by $(f
\rho )g$ and lastly by applying Axiom II. The diagram $\xy
*+{2}*\cir{} \endxy$ commutes by using consecutively the naturality
of $\aaa$, the first diagram of item $(i)$, and both diagram in
equation \eqref{AxIII}. Thus the total diagram commutes, which is sufficient since the left multiplication by $f$ is faithful.
\end{proof}

The converse is not in general true. Under some extra conditions
we can prove

\begin{prop}\label{abelian}
Let $\Sf{B}$ be a locally abelian bicategory whose left and right
multiplications are right exact functors. Consider a WRMC
$\,\Gamma=(f,g,\eta,\rho)$ from $B$ to $A$. If $\eta$ and $\rho$
are epimorphisms, then they are isomorphisms. In particular if
$\eta$ and $\rho$ are epimorphisms, then $A$ and $B$ are
equivalent inside $\Sf{B}$.
\end{prop}
\begin{proof}
We only prove that $\eta:fg \to \II_A$ is an isomorphism. Since
${}_A\Sf{B}_A$ is an abelian category, its suffices to check that
$\eta$ has a null kernel i.e. $\Ker{(\eta)}\,=\,0$. From the sequence
$$
\xy *+{0}="p", p+<2cm,0pt>*+{\Ker{(\eta)}}="1",
p+<4cm,0pt>*+{fg}="2", p+<6cm,0pt>*+{\II_A}="3", {"p" \ar@{->}
"1"}, {"1"\ar@{->}^-{\eta^k} "2"}, {"2" \ar@{->}^-{\eta} "3"}
\endxy
$$
we obtain a commutative diagram
$$
\xy *+{(fg)\Ker{(\eta)}}="p", p+<3cm,0pt>*+{(fg)(fg)}="1",
p+<6cm,0pt>*+{(fg)\II_A}="3", p+<6cm,-2.5cm>*+{\II_A\II_A}="5",
p+<0pt, -2.5cm>*+{\II_A\Ker{(\eta)}}="2",
p+<3cm,-2.5cm>*+{\II_A(fg)}="1'",
p+<2.5cm,-4cm>*+{\Ker{(\eta)}}="4", p+<5.5cm,-4cm>*+{(fg)}="6",
p+<8.5cm,-4cm>*+{\II_A}="8", {"p" \ar@{->}^-{1 \eta^k} "1"}, {"p"
\ar@{->}_-{\eta 1} "2"}, {"1" \ar@{->}^-{1 \eta} "3"}, {"1"
\ar@{->}^-{\eta 1} "1'"}, {"1'" \ar@{->}^-{\lll} "6"} ,{"3"
\ar@{->}^-{\eta 1} "5"}, {"5" \ar@{->}^-{\lll} "8"}, {"2"
\ar@{->}^-{\lll} "4"}, {"4" \ar@{->}_-{\eta^k} "6"}, {"6"
\ar@{->}_-{\eta} "8"}, {"2" \ar@{->}^-{1 \eta^k} "1'"}, {"1'"
\ar@{->}^-{1 \eta} "5"}
\endxy
$$
which implies that $\eta^k \circ \lll_{\Ker{(\eta)}} \circ
(\eta\Ker{(\eta)})\,=\,0$, and so $\eta\Ker{(\eta)}\,=\,0$. Whence
$\Ker{(\eta)}\,=\,0$ since $\eta$ is an epimorphisms and the
multiplication is by hypothesis a right exact functor.
\end{proof}

Let $\Gamma'=(f',g',\eta',\rho')$ be another WRMC from $B$ to $A$.
A \emph{morphism} from $\Gamma$ to $\Gamma'$ is a pair
$(\alpha,\beta): \Gamma \to \Gamma'$ consisting of two $2$-cells
$\alpha: f \to f'$ and $\beta: g \to g'$ rendering commutative the
following diagrams
\begin{equation}\label{W-Morph}
\xy *+{fg}="p", p+<0pt, -2cm>*+{f'g'}="1", p+<2cm,
-1cm>*+{\II_A}="2", {"p" \ar@{->}_-{\alpha\beta} "1"}, {"p"
\ar@{->}^-{\eta} "2"}, {"1" \ar@{->}_-{\eta'} "2"}
\endxy \qquad \qquad \xy *+{gf}="p", p+<0pt, -2cm>*+{g'f'}="1",
p+<2cm, -1cm>*+{\II_B}="2", {"p" \ar@{->}_-{\beta\alpha} "1"},
{"p" \ar@{->}^-{\rho} "2"}, {"1" \ar@{->}_-{\rho'} "2"}
\endxy
\end{equation}
called \emph{compatibility conditions}. The identity morphism is given by the pair of identities $2$-cells
$(1_f,1_g): \Gamma \to \Gamma$. The composition is componentwise. We
thus arrive to
\begin{lemma}\label{lema-1}
The wide right Morita contexts from $B$ to $A$ and their
morphisms form a category which we denote by
${}_A\Sf{W}(\Sf{B})_B$.
\end{lemma}
\begin{proof}
We only need to check that the composition law is well defined.
That is, if $(\alpha,\beta): \Gamma_1 \to \Gamma_2$ and
$(\alpha',\beta'): \Gamma_2 \to \Gamma_3$ are morphisms of WRMC,
then $(\alpha' \circ \alpha,\beta' \circ \beta): \Gamma_1 \to
\Gamma_3$ is also a morphism of WRMC. Thus we need to show that the diagrams of equation \eqref{W-Morph} corresponding to $(\alpha' \circ \alpha,\beta' \circ \beta)$ are commutative. This is clearly
fulfilled as the following commutative diagrams shown
$$
\xy *+{f_1g_1}="p", p+<0pt, -2cm>*+{f_2g_2}="1",
p+<2cm,-2cm>*+{\II_A}="2", p+<0pt,-4cm>*+{f_3g_3}="3", {"p"
\ar@{->}_-{\alpha\beta} "1"}, {"p" \ar@{->}^-{\eta_1} "2"}, {"1"
\ar@{->}_-{\eta_2} "2"}, {"3" \ar@{->}_-{\eta_3} "2"}, {"1"
\ar@{->}_-{\alpha'\beta'} "3"}, {"p"
\ar@/_4pc/|-{(\alpha'\circ\alpha)(\beta'\circ\beta)} "3"}
\endxy \qquad \qquad \xy *+{g_1f_1}="p", p+<0pt, -2cm>*+{g_2f_2}="1",
p+<2cm,-2cm>*+{\II_B}="2", p+<0pt,-4cm>*+{g_3f_3}="3", {"p"
\ar@{->}_-{\beta\alpha} "1"}, {"p" \ar@{->}^-{\rho_1} "2"}, {"1"
\ar@{->}_-{\rho_2} "2"}, {"3" \ar@{->}_-{\rho_3} "2"}, {"1"
\ar@{->}_-{\beta'\alpha'} "3"}, {"p"
\ar@/_4pc/|-{(\beta'\circ\beta)(\alpha'\circ\alpha)} "3"}
\endxy
$$
\end{proof}

Next, we prove that wide right Morita contexts are $1$-cells in a suitable
bicategory. We start by constructing the two sided multiplications
and identities functors.

\subsection*{The functors $\bara{\ccc}_{ABC}$ and
$\bara{\II}_A$} \noindent

Let $A$, $B$ and $C$ are $0$-cells of $\Sf{B}$, and consider the
following WRMC, $\Gamma_i=(f_i,g_i,\eta_i,\rho_i) \in
{}_A\Sf{W}(\Sf{B})_B$ and $\Lambda_i=(p_i,q_i,\gamma_i,\mu_i) \in
{}_B\Sf{W}(\Sf{B})_C$, where $i=1,2$. The $2$-cells are $\eta:fg \to
\II_A$, $\rho:gf \to \II_B$, and $\gamma: pq \to \II_B$, $\mu: qp
\to \II_C$. From these data, we can construct new $2$-cells,
respectively, in the category ${}_A\Sf{B}_A$ and  ${}_C\Sf{B}_C$:
\begin{equation}\label{ast}
\xy *+{(fp)(qg)}="p", p+<2.5cm,0pt>*+{f(p(qg))}="1", p+<5cm,
0pt>*+{f((pq)g)}="2", p+<5cm,-1.5cm>*+{f(\II_Bg)}="3",
p+<5cm,-3cm>*+{fg}="4", p+<2.5cm,-3cm>*+{\II_A}="5", {"p"
\ar@{->}^-{\aaa} "1"}, {"1" \ar@{->}^-{1\,\aaa^{-1}} "2"}, {"2"
\ar@{->}^-{1(\gamma 1)} "3"}, {"3" \ar@{->}^-{1\, \lll} "4"}, {"4"
\ar@{->}_-{\eta} "5"}, {"p" \ar@{-->}_-{\eta \ast \gamma} "5"}
\endxy
\qquad
\xy *+{(qg)(fp)}="p",
p+<2.5cm,0pt>*+{q(g(fp))}="1", p+<5cm, 0pt>*+{q((gf)p)}="2",
p+<5cm,-1.5cm>*+{q(\II_Bp)}="3", p+<5cm,-3cm>*+{qp}="4",
p+<2.5cm,-3cm>*+{\II_C}="5", {"p" \ar@{->}^-{\aaa} "1"}, {"1"
\ar@{->}^-{1\,\aaa^{-1}} "2"}, {"2" \ar@{->}^-{1(\rho 1)} "3"},
{"3" \ar@{->}^-{1\, \lll} "4"}, {"4" \ar@{->}_-{\mu} "5"}, {"p"
\ar@{-->}_-{\mu \ast \rho} "5"}
\endxy
\end{equation}

\begin{lemma}\label{lema-2}
The four-tuples of the form $\Gamma \Lambda:=(fp,qg,\eta \ast \gamma,\mu \ast
\varrho)$ are objects of the category ${}_A\Sf{W}(\Sf{B})_C$.
\end{lemma}
\begin{proof}
We need to verify the commutativity of diagrams in equation
\eqref{Wide} corresponding to $\eta \ast \gamma$ and $\mu \ast
\varrho$. We only check the first one the other is similarly
derived. This diagram decomposes as
\newpage
$$
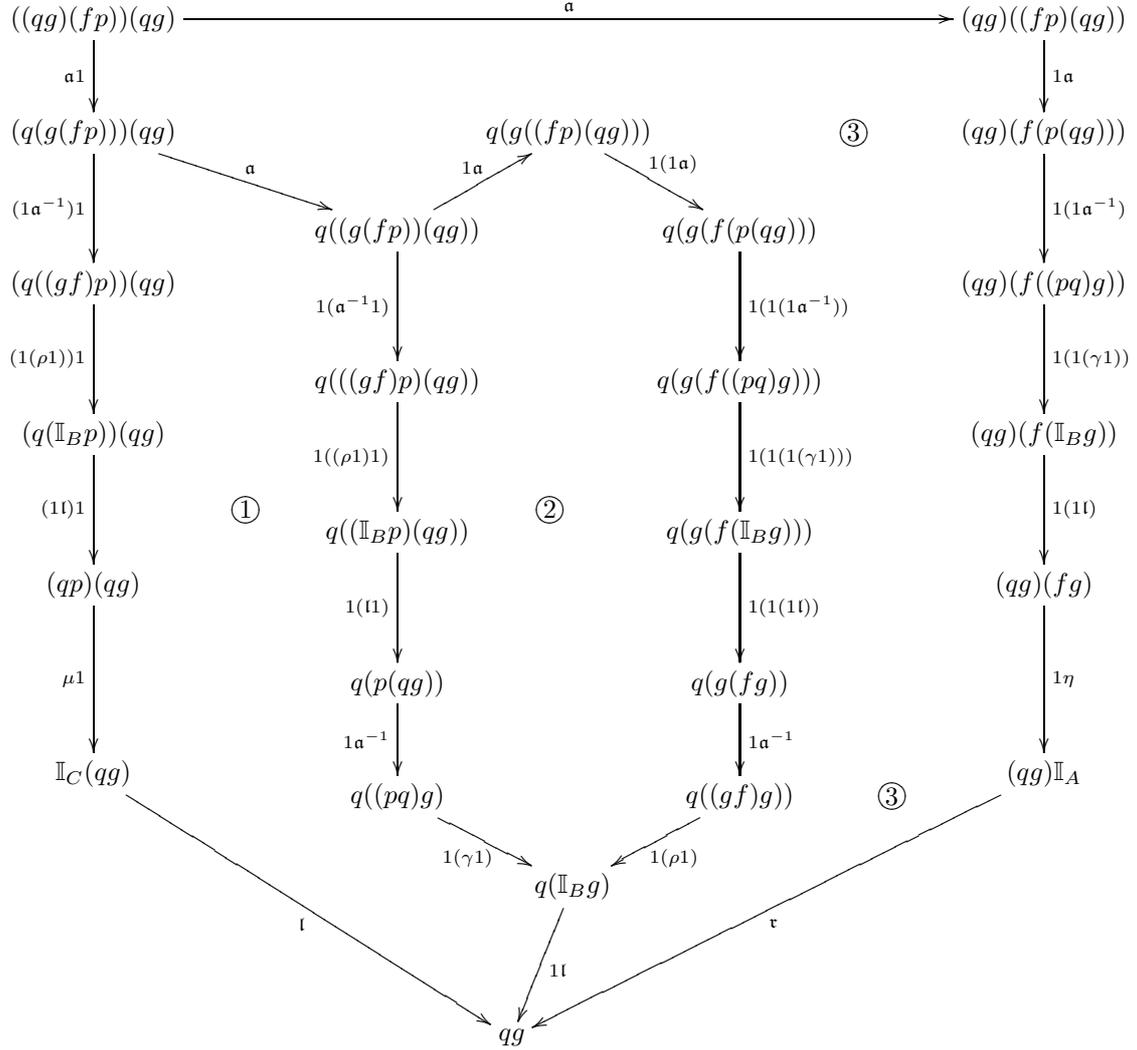
\xy *+{((qg)(fp))(qg)}="p",
p+<12.5cm,0pt>*+{(qg)((fp)(qg))}="1",p+<0pt,-1.5cm>*+{(q(g(fp)))(qg)}="2",
p+<6.25cm,-1.5cm>*+{q(g((fp)(qg)))}="1'",p+<4cm,-2.8cm>*+{q((g(fp))(qg))}="2'",
p+<8.5cm,-2.8cm>*+{q(g(f(p(qg)))}="3'",
p+<12.5cm,-1.5cm>*+{(qg)(f(p(qg)))}="3",
p+<0pt,-3.5cm>*+{(q((gf)p))(qg)}="4",
p+<4cm,-4.8cm>*+{q(((gf)p)(qg))}="4'",
p+<8.5cm,-4.8cm>*+{q(g(f((pq)g)))}="5'",
p+<12.5cm,-3.5cm>*+{(qg)(f((pq)g))}="5",p+<0pt,-5.5cm>*+{(q(\II_Bp))(qg)}="6",
p+<12.5cm,-5.5cm>*+{(qg)(f(\II_Bg))}="7",p+<4cm,-6.8cm>*+{q((\II_B
p)(qg))}="6'",p+<8.5cm,-6.8cm>*+{q(g(f(\II_Bg)))}="7'",
p+<0pt,-7.5cm>*+{(qp)(qg)}="8",p+<12.5cm,-7.5cm>*+{(qg)(fg)}="9",
p+<4cm,-8.8cm>*+{q(p(qg))}="8'", p+<8.5cm, -8.8cm>*+{q(g(fg))}="9'",
p+<0pt,-10cm>*+{\II_C(qg)}="10", p+<4cm,-10.3cm>*+{q((pq)g)}="10'",
p+<12.5cm,-10cm>*+{(qg)\II_A}="11",
p+<8.5cm,-10.3cm>*+{q((gf)g))}="11'",p+<6.3cm,-11.5cm>*+{q(\II_Bg)}="12'",
p+<5.5cm,-13.5cm>*+{qg}="12" ,{"p" \ar@{->}^-{\aaa} "1"},{"p"
\ar@{->}_-{\aaa 1}"2"}, ,{"2" \ar@{->}_-{(1\aaa^{-1})1}"4"},{"4"
\ar@{->}_-{(1(\rho 1))1}"6"},{"6" \ar@{->}_-{(1 \lll)1}"8"},{"8"
\ar@{->}_-{\mu 1}"10"},{"10" \ar@{->}_-{\lll}"12"}, {"1"
\ar@{->}^-{1\aaa} "3"}, {"3" \ar@{->}^-{1(1\aaa^{-1})}"5"},{"5"
\ar@{->}^-{1(1( \gamma1))}"7"},{"7" \ar@{->}^-{1(1\lll)}"9"},{"9"
\ar@{->}^-{1 \eta}"11"},{"11" \ar@{->}^-{\rrr}"12"}, {"2"
\ar@{->}^-{\aaa} "2'"}, {"2'" \ar@{->}^-{1\aaa}"1'"}, {"2'"
\ar@{->}_-{1(\aaa^{-1}1)}"4'"},{"4'" \ar@{->}_-{1((\rho1)1)}"6'"},
{"6'" \ar@{->}_-{1(\lll 1)}"8'"},{"8'" \ar@{->}_-{1\aaa^{-1}}"10'"},
{"10'" \ar@{->}_-{1(\gamma 1)}"12'"}, {"1'"
\ar@{->}^-{1(1\aaa)}"3'"}, {"3'" \ar@{->}^-{1(1(1\aaa^{-1}))}"5'"},
{"5'" \ar@{->}^-{1(1(1(\gamma1)))}"7'"}, {"7'"
\ar@{->}^-{1(1(1\lll))}"9'"}, {"9'" \ar@{->}^-{1\aaa^{-1}}"11'"},
{"11'" \ar@{->}^-{1(\rho1)}"12'"}, {"12'" \ar@{->}^-{1\lll}"12"},
p+<2cm,-6.5cm>*+{1}*\cir{ }, p+<6cm,-6.5cm>*+{2}*\cir{
},p+<10cm,-1.5cm>*+{3}*\cir{ }, p+<10.5cm,-10.3cm>*+{3}*\cir{ }
\endxy
$$

Diagram $\xy *+{1}*\cir{}\endxy$ commutes by applying five times the
naturality of $\aaa$, by Axioms I , II, and the second diagram of
\eqref{AxIII}, and also by the first diagram of equation
\eqref{Wide} satisfied by $\gamma$ and $\mu$. To check the
commutativity of $\xy *+{2}*\cir{}\endxy$, we start by Axiom I, next
apply the naturality of $\aaa$ (in five occasions), and use the
second diagram of \eqref{AxIII} together with the usual obvious fact about
multiplying $2$-cells i.e. $\alpha\beta= 1\beta \circ \alpha
1=\alpha 1 \circ 1\beta$ (used three times). In diagram $\xy
*+{3}*\cir{}\endxy$, we first apply the Axiom I, and continuously
the naturality of $\aaa$ (five times), and ended by Axiom II and the
first diagram of \eqref{Wide} satisfied by $\eta$ and $\rho$. The
total diagram is then commutative and corresponds to the first
diagram of \eqref{Wide} taking the $2$-cells $\eta \ast \gamma$ and
$\mu \ast \rho$.
\end{proof}
Given two morphisms $(\alpha,\beta): \Gamma_1 \to \Gamma_2$ and
$(\varsigma,\tau): \Lambda_1 \to \Lambda_2$, we clearly have two
$2$-cells, namely, $\alpha\varsigma: f_1p_1\to f_2p_2$ and
$\tau\beta:q_1g_1 \to q_2g_2$, which satisfy
\begin{lemma}\label{lema-3}
The pair of $2$-cells $(\alpha\varsigma,\tau\beta)$ defines a
morphism $(\alpha\varsigma,\tau\beta): \Gamma_1\Lambda_1 \to
\Gamma_2\Lambda_2$ in the category ${}_A\Sf{W}(\Sf{B})_C$, where
$\Gamma_i \Lambda_i$ are defined as in Lemma \ref{lema-2},
$i=1,2$.
\end{lemma}
\begin{proof}
The first diagram in equation \eqref{W-Morph} corresponding to the stated
pair of morphisms commutes as it is shown by the following diagram
$$
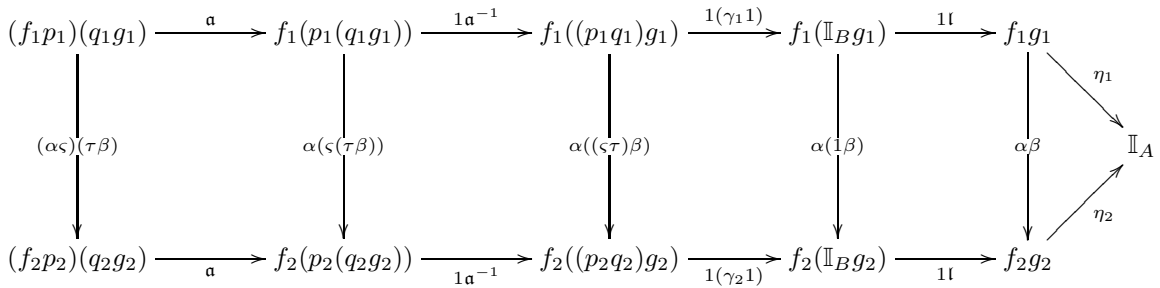
\xy *+{(f_1p_1)(q_1g_1)}="p", p+<3.5cm,
0pt>*+{f_1(p_1(q_1g_1))}="1", p+<7cm,
0pt>*+{f_1((p_1q_1)g_1)}="3", p+<10cm, 0pt>*+{f_1(\II_Bg_1)}="5",
p+<12.5cm, 0pt>*+{f_1g_1}="7", p+<14cm, -1.5cm>*+{\II_A}="9",
p+<0pt,-3cm>*+{(f_2p_2)(q_2g_2)}="2",
p+<3.5cm,-3cm>*+{f_2(p_2(q_2g_2))}="4",
p+<7cm,-3cm>*+{f_2((p_2q_2)g_2)}="6",
p+<10cm,-3cm>*+{f_2(\II_Bg_2)}="8",
p+<12.5cm,-3cm>*+{f_2g_2}="10", {"p" \ar@{->}^-{\aaa} "1"}, {"1"
\ar@{->}^-{1\aaa^{-1}} "3"}, {"1"
\ar@{->}|-{\alpha(\varsigma(\tau\beta))} "4"}, {"3"
\ar@{->}^-{1(\gamma_11)} "5"}, {"3"
\ar@{->}|-{\alpha((\varsigma\tau)\beta)} "6"}, {"5"
\ar@{->}^-{1\lll} "7"}, {"5" \ar@{->}|-{\alpha(1\beta)} "8"}, {"7"
\ar@{->}^-{\eta_1} "9"}, {"7" \ar@{->}|-{\alpha\beta} "10"}, {"p"
\ar@{->}|-{(\alpha\varsigma)(\tau\beta)} "2"},{"2"
\ar@{->}_-{\aaa} "4"}, {"4" \ar@{->}_-{1\aaa^{-1}} "6"}, {"6"
\ar@{->}_-{1(\gamma_21)} "8"}, {"8" \ar@{->}_-{1\lll} "10"}, {"10"
\ar@{->}_-{\eta_2} "9"}
\endxy
$$
The commutativity of the second one is similarly obtained.
\end{proof}
Since the multiplication $\ccc$ of two identities $2$-cells gives an
identity, the above multiplication respects also this rule. By
lemmata \ref{lema-2} and \ref{lema-3}, we thus get the following
\begin{cor}\label{baraC}
Given $A$, $B$ and $C$ three $0$-cells of $\Sf{B}$, there is a
covariant functor
$$
\xymatrix@R=0pt{\bara{\ccc}_{ABC}:\,\, {}_A\Sf{W}(\Sf{B})_B \times
{}_B\Sf{W}(\Sf{B})_C \ar@{->}[rr] & & {}_A\Sf{W(B)}_C \\
\lr{(\Gamma, \Lambda)} \ar@{->}[rr] & &
\Gamma\Lambda\,=\,\lr{fp,qg,\eta\ast\gamma,\mu\ast\varrho} \\
\lr{(\alpha,\beta), (\varsigma,\tau)} \ar@{->}[rr] & &
(\alpha,\beta)(\varsigma,\tau) \,=\, (\alpha\varsigma,\tau\beta)}
$$
where $\Gamma=(f,g,\eta,\rho)$ and $\Lambda=(p,q,\gamma,\mu)$.
\end{cor}

The identity functors are given as follows. Let $A$ be any
$0$-cell of $\Sf{B}$, define the four-tuple
$\bara{\II}_A=(\II_A,\II_A, \rrr_{\II_A}, \lll_{\II_A})$, then we
have the following diagrams
$$
\xy *+{\II_A(\II_A\II_A)}="p", p+<2.5cm,0pt>*+{\II_A\II_A}="1",
p+<4cm,-1cm>*+{\II_A}="2", p+<0pt,-2cm>*+{(\II_A\II_A)\II_A}="3",
p+<2.5cm,-2cm>*+{\II_A \II_A}="4", {"p" \ar@{->}^-{1 \rrr} "1"},
{"1" \ar@{->}^-{\rrr } "2"}, {"3" \ar@{->}^-{\aaa} "p"}, {"3"
\ar@{->}_-{\lll 1} "4"}, {"4" \ar@{->}_-{\lll} "2"}
\endxy \qquad \qquad \xy *+{(\II_A\II_A)\II_A}="p",
p+<2.5cm,0pt>*+{\II_A \II_A}="1", p+<4cm,-1cm>*+{\II_A}="2",
p+<0pt,-2cm>*+{\II_A(\II_A\II_A)}="3", p+<2.5cm,-2cm>*+{\II_A\II_A
}="4" , {"p" \ar@{->}^-{\rrr 1} "1"}, {"1" \ar@{->}^-{\lll} "2"},
{"p" \ar@{->}_-{\aaa} "3"}, {"3" \ar@{->}_-{1 \lll} "4"}, {"4"
\ar@{->}_-{\rrr} "2"}
\endxy
$$ where the second diagram commutes by Axiom II, and the first one decomposes as
$$
\xy *+{\II_A(\II_A\II_A)}="p", p+<6.5cm,0pt>*+{\II_A\II_A}="1",
p+<8cm,-1.5cm>*+{\II_A}="2", p+<0pt,-3cm>*+{(\II_A\II_A)\II_A}="3",
p+<6.5cm,-3cm>*+{\II_A \II_A}="4", p+<4cm,-1.5cm>*+{\II_A
\II_A}="2'", {"p" \ar@{->}^-{1 \rrr} "1"}, {"1" \ar@{->}^-{\rrr }
"2"}, {"3" \ar@{->}^-{\aaa} "p"}, {"3" \ar@{->}_-{\lll 1} "4"}, {"4"
\ar@{->}_-{\lll} "2"}, {"2'" \ar@{->}_-{\rrr} "2"}, {"p"
\ar@{->}^-{\rrr} "2'"}, {"3" \ar@{->}_-{\rrr 1} "2'"}
\endxy
$$
which is commutative by the first diagram of \eqref{AxIII}, the
natutrality of $\rrr$, and the equality
$\rrr_{\II_A}=\lll_{\II_A}$.

Now, we look at the associativity of the two sided multiplications
$\bara{\ccc}$.

\subsection*{The natural isomorphisms $\bara{\aaa}_{\Gamma\Lambda\Omega}$}
\noindent

Let $\Gamma=(f,g,\eta,\varrho) \in {}_A\Sf{W}(\Sf{B})_B$,
$\Lambda=(p,q,\gamma,\mu) \in {}_B\Sf{W}(\Sf{B})_C$, and
$\Omega=(u,v,\theta,\sigma) \in {}_C\Sf{W}(\Sf{B})_D$. We then
have the following $2$-cells:
$$
\xymatrix@R=0pt{\eta:\,fg \ar@{->}[r] & \II_A \\ \rho:\,gf
\ar@{->}[r] & \II_B, }\qquad \xymatrix@R=0pt{\gamma:\,
pq\ar@{->}[r] & \II_B \\ \mu:\,qp \ar@{->}[r] & \II_C, }\qquad
\xymatrix@R=0pt{\theta:\,uv \ar@{->}[r] & \II_C \\ \sigma:\,vu
\ar@{->}[r] & \II_D, }
$$ where each vertical pair satisfies \eqref{Wide}.
By Corollary \ref{baraC}, we have two WRMC
$$
\lr{\Gamma \Lambda} \Omega \,=\, (fp,qg,\eta\ast\gamma,\mu\ast\rho)(u,v,\theta,\sigma)
\,=\, \lr{(fp)u, v(qg), (\eta\ast\gamma)\ast \theta, \sigma\ast(\mu\ast\rho)} \in {}_A\Sf{W}(\Sf{B})_D
$$
and
$$
\Gamma \lr{\Lambda \Omega} \,=\, (f,g,\eta,\rho)(pu, vq, \gamma \ast\theta,\sigma \ast \mu)
\,=\, (f(pu), (vq)g, \eta \ast (\gamma \ast\theta),(\sigma \ast  \mu)\ast \rho) \in {}_A\Sf{W}(\Sf{B})_D
$$

\begin{lemma}
The pair
$\bara{\aaa}_{\Gamma\Lambda\Omega}:=(\aaa,\aaa^{-1}):\lr{\Gamma
\Lambda} \Omega \to \Gamma \lr{\Lambda \Omega}$ defines an
isomorphism in the category ${}_A\Sf{W}(\Sf{B})_D$ with inverse
$(\aaa^{-1},\aaa)$. Moreover, $\bara{\aaa}_{---}: (- -)- \longrightarrow -(- -)$ is a natural
isomorphism.
\end{lemma}
\begin{proof}
For the first statement, we need to show that the diagrams in
equation \eqref{W-Morph} associated to the pair $(\aaa,\aaa^{-1})$
are commutative. We only prove the commutativity of the first one
the remainder is similarly deduced. This first diagram is
$$\xy
*+{((fp)u)(v(qg))}="p", p+<6cm,0pt>*+{(f(pu))((vq)g)}="1",
p+<3cm,-2cm>*+{\II_A}="2", {"p" \ar@{->}^-{\aaa\aaa^{-1}} "1"},
{"p" \ar@{->}_-{(\eta \ast\gamma)\ast \theta} "2"}, {"1"
\ar@{->}^-{\eta \ast(\gamma\ast \theta)} "2"}
\endxy
$$ writing explicitly the maps, we get
\newpage

$$
\xy *+{((fp)u)(v(qg))}="p", p+<13cm,0pt>*+{(f(pu))((vq)g)}="1",
p+<0pt,-2.5cm>*+{(fp)(u(v(qg)))}="2",
p+<0pt,-5cm>*+{(fp)((uv)(qg))}="4",
p+<4cm,-6.3cm>*+{(fp)(((uv)q)g)}="4'",
p+<0pt,-7.5cm>*+{(fp)(\II_C(qg))}="6",
p+<0pt,-9.5cm>*+{(fp)(qg)}="8", p+<0pt,-11.5cm>*+{f(p(qg))}="10",
p+<6cm,-12.5cm>*+{f((pq)g)}="12",
p+<2.5cm,-13.8cm>*+{f(\II_Bg)}="13", p+<6cm,-13.8cm>*+{fg}="14",
p+<8cm,-13.8cm>*+{\II_A}="15", {"12" \ar@{->}|-{1(\gamma 1)}
"13"}, {"13" \ar@{->}^-{1\lll} "14"}, {"14" \ar@{->}^-{\eta}
"15"}, p+<13cm,-2.5cm>*+{f((pu)((vq)g))}="3",
p+<13cm,-5cm>*+{f(((pu)(vq))g)}="5",
p+<13cm,-7.5cm>*+{f((p(u(vq)))g)}="7",
p+<13cm,-9.5cm>*+{f((p((uv)q))g)}="9",
p+<13cm,-11.5cm>*+{f((p(\II_Cq))g)}="11",
p+<8cm,-2.5cm>*+{f(p(u((vq)g)))}="3'",
p+<8cm,-5cm>*+{f(p((u(vq))g))}="5'",
p+<8cm,-7.5cm>*+{f(p(((uv)q)g))}="7'",
p+<8cm,-9.5cm>*+{f(p((\II_Cq)g))}="9'",
{"p"\ar@{->}^-{\aaa\aaa^{-1}} "1"}, {"p" \ar@{->}|-{\aaa} "2"},
{"2" \ar@{->}|-{1\aaa^{-1}} "4"}, {"4" \ar@{->}|-{1(\theta 1)}
"6"}, {"6" \ar@{->}|-{1\lll} "8"}, {"8" \ar@{->}|-{\aaa} "10"},
{"10" \ar@{->}|-{1\aaa^{-1}} "12"}, {"1" \ar@{->}|-{\aaa} "3"},
{"3" \ar@{->}|-{1\aaa^{-1}} "5"}, {"5" \ar@{->}|-{1(\aaa1)} "7"},
{"7" \ar@{->}|-{1((1\aaa^{-1})1)} "9"}, {"9"
\ar@{->}|-{1((1(\theta 1))1)} "11"}, {"11" \ar@{->}|-{1((1\lll)1)}
"12"}, {"3" \ar@{->}|-{1\aaa} "3'"}, {"3'"
\ar@{->}|-{1(1\aaa^{-1})} "5'"}, {"5'"
\ar@{->}|-{1(1(\aaa^{-1}1))} "7'"}, {"7'" \ar@{->}|-{1(1((\theta
1)1))} "9'"}, {"9'" \ar@{->}|-{1(1(\lll 1))} "10"}, {"4"
\ar@{->}|-{1\aaa^{-1}} "4'"}, {"4'" \ar@{->}|-{\aaa} "7'"},
p+<4cm,-2cm>*+{1}*\cir{ }, p+<4cm,-8.5cm>*+{2}*\cir{ },
p+<10.5cm,-10cm>*+{3}*\cir{ }
\endxy
$$
Diagram $\xy *+{1}*\cir{}\endxy$ commutes by Axiom I (applied two
times) and by naturality of $\aaa$ (used three times). Applying
the second diagram of \eqref{AxIII}, and three times the
naturality of $\aaa$, we get the commutativity of $\xy
*+{2}*\cir{}\endxy$. It is clear that $\xy *+{3}*\cir{}\endxy = f
\times \xy *+{3'}*\cir{}\endxy$, and $\xy *+{3'}*\cir{}\endxy$ is
commutative before applying Axiom I, and consecutively three
times the naturality of $\aaa$.\\ The rest of the proof is clearly
obtained by using the natural isomorphism $\aaa$.
\end{proof}

\subsection*{The natural isomorphisms $\bara{\rrr}_{\Gamma}$ and
$\bara{\lll}_{\Gamma}$} \noindent

Let $A$ and $B$ two $0$-cells of $\Sf{B}$, we have shown that
$\bara{\II}_A=(\II_A,\II_A, \rrr_{\II_A},\lll_{\II_A})$ and
$\bara{\II}_B=(\II_B,\II_B, \rrr_{\II_B},\lll_{\II_B})$ are,
respectively, objects of ${}_A\Sf{W}(\Sf{B})_A$ and
${}_B\Sf{W}(\Sf{B})_B$. Consider now $\Gamma=(f,g,\eta,\rho)$
 any object of  the category ${}_A\Sf{W}(\Sf{B})_B$, using the multiplication $\bara{\ccc}$,
we can define two $2$-cells: $$(\rrr_f, \lll_g): \,\Gamma
\bara{\II}_B=(f\II_B,\II_Bg, \eta \ast \rrr, \lll \ast \rho)
\longrightarrow \Gamma, \quad (\lll_g,\rrr_f): \,\bara{\II}_A
\Gamma=(\II_Af,g\II_A,\rrr \ast \eta, \rho \ast \lll)
\longrightarrow \Gamma$$
\begin{lemma}
Keep the above notation.
\begin{enumerate}[(i)]
\item The pairs $(\rrr_f, \lll_g):\,\Gamma \bara{\II}_B \to
\Gamma$ and $(\lll_g,\rrr_f):\, \bara{\II}_A \Gamma \to \Gamma$
are isomorphisms of the category ${}_A\Sf{W}(\Sf{B})_B$ with
inverses, respectively, $(\rrr_f^{-1}, \lll_g^{-1})$ and
$(\lll_g^{-1},\rrr_f^{-1})$.

\item These in fact define natural isomorphisms
$\bara{\rrr}_{-}$ and $\bara{\lll}_{-}$ given, respectively, at the object
$\Gamma=(f,g,\eta,\rho)$ by $\bara{\rrr}_{\Gamma}\,=\,(\rrr_f,
\lll_g)$ and  $\bara{\lll}_{\Gamma}\,=\,(\lll_g,\rrr_f)$.
\end{enumerate}
\end{lemma}
\begin{proof}
$(i)$ We only prove that $(\rrr_f,\lll_g)$ is a morphisms of
${}_A\Sf{W}(\Sf{B})_B$. The proof of the other morphism is
similarly given. We thus need to show that the following diagrams
are commutative
\begin{equation}\label{Prf-1}
 \xy *+{(f\II_B)(\II_Bg)}="p", p+<0pt, -2.5cm>*+{fg}="1", p+<2.5cm,
-1.25cm>*+{\II_A}="2", {"p" \ar@{->}_-{\rrr\lll} "1"}, {"p"
\ar@{->}^-{\eta \ast \rrr_{\II_B}} "2"}, {"1" \ar@{->}_-{\eta}
"2"}
\endxy \qquad \qquad \xy *+{(\II_Bg)(f\II_B)}="p", p+<0pt,
-2.5cm>*+{gf}="1", p+<2.5cm, -1.25cm>*+{\II_B}="2", {"p"
\ar@{->}_-{\lll\rrr} "1"}, {"p" \ar@{->}^-{\lll_{\II_A} \ast \rho}
"2"}, {"1" \ar@{->}_-{\rho} "2"}
\endxy
\end{equation}
The first diagram decomposes as
$$
\xy *+{(f\II_B)(\II_Bg)}="p", p+<3.5cm,
0pt>*+{f(\II_B(\II_Bg))}="1", p+<7cm, 0pt>*+{f((\II_B
\II_B)g)}="2", p+<0pt, -2cm>*+{(f\II_B)g}="3", p+<3.5cm,
-2cm>*+{f(\II_Bg)}="4", p+<7cm, -3.5cm>*+{\II_A}="5", p+<3.5cm,
-3.5cm>*+{fg}="6", {"p" \ar@{->}^-{\aaa} "1"}, {"1"
\ar@{->}^-{1\aaa^{-1}} "2"}, {"p" \ar@{->}_-{1\lll} "3"}, {"3"
\ar@{->}^-{\aaa} "4"}, {"2" \ar@{->}^-{1(\rrr1)} "4"}, {"1"
\ar@{->}^-{1(1\lll)} "4"}, {"3" \ar@{->}_-{\rrr 1} "6"}, {"4"
\ar@{->}^-{1\lll} "6"}, {"6" \ar@{->}_-{\eta} "5"}
\endxy
$$
which is commutative by naturality of $\aaa$, the first diagram in
equation \eqref{AxIII}, and Axiom II.\\ The second diagram in
\eqref{Prf-1} decomposes as
$$
\xy +*{(\II_Bg)(f\II_B)}="p",
p+<3.5cm,0pt>*+{\II_B(g(f\II_B))}="1", p+<7.5cm,
0pt>*+{\II_B((fg)\II_B)}="3", p+<11cm,
0pt>*+{\II_B(\II_B\II_B)}="5", p+<0pt, -2cm>*+{(\II_Bg)f}="2",
p+<2.5cm, -2cm>*+{\II_B(fg))}="1'", p+<0pt, -5.5cm>*+{fg}="4",
p+<12.5cm, -7cm>*+{\II_B\II_B}="4'", p+<8cm, -10cm>*+{\II_B}="6",
p+<6cm, -4cm>*+{(\II_B(gf))\II_B}="3'", p+<9.5cm,
-4cm>*+{(\II_B\II_B)\II_B}="2'", p+<5cm, -6cm>*+{(gf)\II_B}="5'",
p+<8.5cm, -6cm>*+{\II_B\II_B}="6'", {"p" \ar@{->}^-{\aaa} "1"},
{"p" \ar@{->}_-{1 \rrr} "2"}, {"1" \ar@{->}^-{1\aaa^{-1}} "3"},
{"3" \ar@{->}^-{1(\rho 1)} "5"}, {"3" \ar@{->}^-{1\rrr} "1'"},
{"1'" \ar@{->}_-{\lll} "4"}, {"3" \ar@{->}^-{\aaa^{-1}} "3'"},
{"5" \ar@{->}^-{\aaa^{-1}} "2'"}, {"5" \ar@{->}^-{1\lll} "4'"},
{"3'" \ar@{->}^-{\lll 1} "5'"}, {"2'" \ar@{->}^-{\lll 1} "6'"},
{"3'" \ar@{->}^-{(1\rho)1} "2'"}, {"5'" \ar@{->}_-{\rho 1} "6'"},
{"4'" \ar@{->}_-{\lll} "6"}, {"2" \ar@{->}_-{\lll 1} "4"}, {"4"
\ar@{->}_-{\rho} "6"}, {"6'" \ar@{->}^-{\rrr} "6"}, {"5'"
\ar@{->}^-{\rrr} "4"}, p+<2cm, -1cm>*+{1}*\cir{}, p+<3.5cm,
-3.5cm>*+{2}*\cir{}, p+<8.5cm, -2cm>*+{3}*\cir{}, p+<7.75cm,
-5cm>*+{4}*\cir{}, p+<11cm, -5.5cm>*+{5}*\cir{}, p+<6.5cm,
-7.5cm>*+{6}*\cir{}
\endxy
$$
Diagram $\xy *+{1}*\cir{} \endxy$ commutes by naturality of $\aaa$
and by the second diagram in equation \eqref{AxIII}. Applying
equation \eqref{AxIV}, we get the commutativity of $\xy
*+{2}*\cir{}
\endxy$. The diagrams $\xy *+{3}*\cir{} \endxy$, $\xy *+{4}*\cir{}
\endxy$ and $\xy *+{6}*\cir{} \endxy$ are commutative, since
$\aaa^{-1}$, $\lll$ and $\rrr$ are natural transformations.
Lastly, diagram $\xy *+{5}*\cir{} \endxy$ commutes by the second
diagram in \eqref{AxIII} and the equality $\lll_{\II_A} =
\rrr_{\II_A}$.\\ $(ii)$ It is componentwise derived from the
naturaliy of $\rrr_{-}$ and $\lll_{-}$.
\end{proof}

The Axioms $\bara{\rm I}$ and $\bara{ \rm II}$ corresponding to
the functors $\bara{\ccc}$, $\bara{\II}$ and the natural
isomorphisms $\bara{\aaa}$, $\bara{\rrr}$, and $\bara{\lll}$, are
easily derived from Axiom I and II. We then have
\begin{prop}\label{Bicategory}
Let $\Sf{B}$ be a bicategory. Then the following data form a
bicategory $\Sf{W}(\Sf{B})$

\begin{enumerate}[$\bullet$]
\item $0$-cells. Are all $0$-cells of $\Sf{B}$

\item $1$-cells. From $B$ to $A$ are four-tuples
$\Gamma=(f,g,\eta,\rho)$ consisting of two $1$-cells $f \in
{}_A\Sf{B}_B$ and $g \in {}_B\Sf{B}_A$, and two $2$-cells $\eta:
fg \to \II_A$ and $\rho: gf \to \II_B$ satisfying the
compatibility conditions (i.e. the diagrams of \eqref{Wide} are
commutative). The identity $1$-cell is given by the four-tuple
$\bara{\II}_A\,=\,(\II_A,\II_A,\rrr_{\II_A},\lll_{\II_A})$.

\item $2$-cells. Are pairs $(\alpha,\beta): \Gamma=(f,g,\eta,\varrho) \to
\Gamma'=(f',g',\eta',\varrho')$ consisting of two $2$-cells
$\alpha:f \to f'$ and $\beta: g \to g'$ rendering commutative the
diagrams of equation \eqref{W-Morph}.
\end{enumerate}
\end{prop}

The construction of wide right Morita contexts is in fact functorial.

\begin{prop}\label{Morph-bicatg}
Let  $(\Scr{F},\Phi): \Sf{B} \to \Sf{B'}$ be a homomorphism of bicategories. Then there is a homomorphism of bicategories $\Sf{W}({\Scr{F}}, \Phi) : \Sf{W(B)} \to \Sf{W(B')}$.
\end{prop}
\begin{proof}
The function on $0$-cell coincides with that of $\Scr{F}$. Let $\Gamma=(f,g,\eta,\rho)$ be any $1$-cell in ${}_A\Sf{W(B)}_B$. Define the following two $2$-cells in ${}_{\Scr{F}(A)}\Sf{B'}_{\Scr{F}(B)}$
$$\xymatrix@C=40pt{ \eta': \Scr{F}(f)\Scr{F}(g) \ar@{->}^-{\Phi^1_{fg}}[r] & \Scr{F}(fg) \ar@{->}^-{\Scr{F}(\eta)}[r] & \Scr{F}(\II_A)   \ar@{->}^-{\Phi^0{ }^{-1}}[r] & \II'_{\Scr{F}(A)}, }
$$
$$\xymatrix@C=40pt{ \rho': \Scr{F}(g)\Scr{F}(f) \ar@{->}^-{\Phi^1_{gf}}[r] & \Scr{F}(gf) \ar@{->}^-{\Scr{F}(\rho)}[r] & \Scr{F}(\II_B)   \ar@{->}^-{\Phi^0{ }^{-1}}[r] & \II'_{\Scr{F}(B)}. }
$$
Set $\Gamma'=(\Scr{F}(f),\Scr{F}(g),\eta',\rho')$, we claim that $\Gamma' \in {}_{\Scr{F}(A)}\Sf{W(B')}_{\Scr{F}(B)}$. We need to check that the diagrams of equation \eqref{Wide} corresponding to $\eta'$ and $\rho'$ in $\Sf{B'}$ are commutative. The first one decomposes as
$$
\xy
+*{\Scr{F}(g)(\Scr{F}(f)\Scr{F}(g))}="p",
p+<0pt, -2cm>*+{(\Scr{F}(g)\Scr{F}(f))\Scr{F}(g)}="1", p+<0pt, -4cm>*+{\Scr{F}(gf)\Scr{F}(g)}="2", p+<0pt,-6cm>*+{\Scr{F}(\II_B)\Scr{F}(g)}="3", p+<5cm,-8cm>*+{\II'_{\Scr{F}(B)}\Scr{F}(g)}="4", p+<4cm,0pt>*+{\Scr{F}(g)\Scr{F}(fg)}="5", p+<8cm,0pt>*+{\Scr{F}(g)\Scr{F}(\II_A)}="6", p+<12cm,0pt>*+{\Scr{F}(g)\II'_{\Scr{F}(A)}}="7", p+<12cm,-6cm>*+{\Scr{F}(g)}="8", p+<4cm,-2cm>*+{\Scr{F}(g(fg))}="15", p+<8cm,-2cm>*+{\Scr{F}(g\II_{A})}="16", p+<4cm,-4cm>*+{\Scr{F}((gf)g)}="25", p+<8cm,-6cm>*+{\Scr{F}(\II_{B}g)}="36", { "1" \ar@{->}^-{\aaa'} "p"}, { "1" \ar@{->}_-{\Phi 1} "2"}, { "2" \ar@{->}_-{\Scr{F}(\rho) 1} "3"},{ "2" \ar@{->}^-{\Phi} "25"}, { "3" \ar@{->}_-{\Phi^{-1} 1} "4"}, { "4" \ar@{->}_-{\lll'} "8"}, { "p" \ar@{->}^-{1 \Phi} "5"}, { "5" \ar@{->}^-{1\Scr{F}(\eta)} "6"}, { "5" \ar@{->}^-{\Phi} "15"}, { "6" \ar@{->}^-{1 \Phi^{-1}} "7"}, { "6" \ar@{->}_-{\Phi} "16"}, { "7" \ar@{->}^-{\rrr'} "8"}, { "15" \ar@{->}^-{\Scr{F}(1 \eta)} "16"}, { "25" \ar@{->}^-{\Scr{F}(\aaa)} "15"}, { "25" \ar@{->}^-{\Scr{F}(\rho 1)} "36"}, { "3" \ar@{->}^-{\Phi} "36"}, { "16" \ar@{->}^-{\Scr{F}(\rrr)} "8"}, { "36" \ar@{->}^-{\Scr{F}(\lll)} "8"},
\endxy
$$
which is clearly a commutative diagram. Similarly we check that the second one is also commutative.

Now, given a $2$-cell $(\alpha,\beta): \Gamma \to \Sigma=(h,e, \sigma, \tau)$ in ${}_A\Sf{W(B)}_B$, we get by the naturality of $\Phi$ and equation \eqref{W-Morph} satisfied by the pair $(\alpha, \beta)$, commutative diagrams
$$
\xymatrix@C=30pt@R=10pt{ \Scr{F}(f)\Scr{F}(g) \ar@{->}^-{\Phi}[r] \ar@{->}_-{\Scr{F}(\alpha)\Scr{F}(\beta)}[dd]  & \Scr{F}(fg) \ar@{->}^-{\Scr{F}(\eta)}[r] \ar@{->}^-{\Scr{F}(\alpha\beta)}[dd] & \Scr{F}(\II_A) \ar@{->}^-{\Phi^{-1}}[rd] &  \\ & &  & \II'_{\Scr{F}(A)} \\ \Scr{F}(h)\Scr{F}(e) \ar@{->}_-{\Phi}[r] & \Scr{F}(he) \ar@{->}_-{\Scr{F}(\sigma)}[r] & \Scr{F}(\II_A) \ar@{->}_-{\Phi^{-1}}[ru] & }
$$
and
$$
\xymatrix@C=30pt@R=10pt{ \Scr{F}(g)\Scr{F}(f) \ar@{->}^-{\Phi}[r] \ar@{->}_-{\Scr{F}(\beta)\Scr{F}(\alpha)}[dd]  & \Scr{F}(gf) \ar@{->}^-{\Scr{F}(\rho)}[r] \ar@{->}^-{\Scr{F}(\beta\alpha)}[dd]  & \Scr{F}(\II_B) \ar@{->}^-{\Phi^{-1}}[rd] &  \\ & &  & \II'_{\Scr{F}(B)} \\ \Scr{F}(e)\Scr{F}(h) \ar@{->}_-{\Phi}[r] & \Scr{F}(eh) \ar@{->}_-{\Scr{F}(\tau)}[r] & \Scr{F}(\II_B) \ar@{->}_-{\Phi^{-1}}[ru] & }
$$
This means that $(\Scr{F}(\alpha),\Scr{F}(\beta)): (\Scr{F}(f),\Scr{F}(g),\eta',\rho') \to (\Scr{F}(h),\Scr{F}(e), \sigma', \tau')$ is a $2$-cell in ${}_{\Scr{F}(A)} \Sf{W(B')}_{\Scr{F}(B)}$. We have thus construct a family of functors
$$
\xymatrix@C=50pt@R=0pt{ {}_A\Sf{W(B)}_B \ar@{->}^-{\Sf{W}(\Scr{F})_{A,\,B}}[rr] & & {}_{\Scr{F}(A)}\Sf{W(B')}_{\Scr{F}(B)} \\ \Gamma=(f,g,\eta,\rho) \ar@{->}[rr] & & \Sf{W}(\Scr{F})(\Gamma):=(\Scr{F}(f),\Scr{F}(g),\eta',\rho') \\ \left[ \underset{}{} (\alpha,\beta): \Gamma \to \Sigma \right] \ar@{->}[rr] & & \left[ \underset{}{}  \Sf{W}(\Scr{F})(\alpha,\beta):=(\Scr{F}(\alpha),\Scr{F}(\beta)): \Sf{W}(\Scr{F})(\Gamma) \to \Sf{W}(\Scr{F})(\Sigma) \right]}
$$
On the other hand, given $\Gamma=(f,g,\eta,\rho)$ and $\Lambda=(p,q,\mu,\gamma)$ two $1$-cells, respectively, in ${}_A\Sf{W(B)}_B$ and ${}_B\Sf{W(B)}_C$, using the naturality of $\Phi$, we obtain two $2$-cells
$$
\xymatrix@R=4pt@C=50pt{\Sf{W}(\Scr{F})(\Gamma)\, \Sf{W}(\Scr{F})(\Lambda) \ar@{=}[d] & & \Sf{W}(\Scr{F})(\Gamma \Lambda) \ar@{=}[d] \\
\lr{\Scr{F}(f)\Scr{F}(p), \Scr{F}(q) \Scr{F}(g), \eta' * \mu', \rho'*\gamma'} \ar@{->}^-{(\Phi_{fp}{}^{-1}, \Phi_{qg}{}^{-1})}[rr] & & \lr{\Scr{F}(fp), \Scr{F}(qg), (\eta*\mu)', (\rho*\gamma)'} }
$$
and
$$
\xymatrix@R=3pt@C=50pt{\Sf{W}(\Scr{F})\lr{\II_A,\II_A,\rrr_A,\lll_A} \ar@{=}[d] & & \lr{\II'_{\Scr{F}(A)},\II'_{\Scr{F}(A)},\rrr'_{\II'_{\Scr{F}(A)}},\lll'_{\II'_{\Scr{F}(A)}} } \\
\lr{\Scr{F}(\II_A), \Scr{F}(\II_A), \rrr'_{\II_A}, \lll'_{\II_A}} \ar@{->}_-{(\Phi_{A}{}^{-1}, \Phi_{A}{}^{-1})}[rru] & &  }
$$
The axioms satisfied by the natural transformations $\Phi$ show that $\lr{\Sf{W}(\Scr{F}), (\Phi^{-1},\Phi^{-1})}: \Sf{W(B)} \to \Sf{W(B')}$ is a homomorphism of bicategories as claimed.
\end{proof}

\section{Application to $\Sf{REM(Bim)}$}

In what follows all rings are assumed to be associative with $1$.
Modules are unital modules, and bimodules are left and right
unital modules. The identity linear map associated to any module
$X$ is denoted by the module itself $X$. Given $A$ and $B$ two
rings, the category of $(A, B)$-bimodules (left $A$-modules and
right $B$-modules) is denoted as usual by $\Bimod{A}{B}$. The
symbol $-\tensor{A}-$ between bimodules and bilinear map denotes
the tensor product over $A$. The bimodules bicategory $\Sf{Bim}$
has the class of $0$-cells all rings $A$, $B$, $C$..., and
Hom-categories the categories of bimodules $\Bimod{A}{B}$. The
vertical and horizontal multiplications are given by the tensor
product over rings. For every $M \in \Bimod{A}{B}$, we denote by
$\iota_{M}$ the obvious natural isomorphism at $M$, $\iota_M:
A\tensor{A}M \to M \tensor{B}B$, that is, $\iota_M={\bf r}_M^{-1}
\circ {\bf l}_M$, where ${\bf r}$ and ${\bf l}$ are respectively,
the right and left identities natural isomorphism of $\Sf{Bim}$.

Let $A$ be a ring, an \emph{$A$-coring} \cite{Sweedler:1975} is a
three-tuple $(\coring{C}, \Delta, \varepsilon)$ consisting of an
$A$-bimodule and two $A$-bilinear maps $$\Delta: \coring{C} \to
\coring{C} \tensor{A} \coring{C}\quad \textrm{and}\quad
\varepsilon: \coring{C} \to A,$$ known as the comultiplication and
the counit of $\cc$, which satisfy
$$(\coring{C}\tensor{A}\Delta) \circ \Delta \,\, =\,\, (\Delta
\tensor{A} \coring{C}) \circ \Delta,\quad (\coring{C} \tensor{A}
\varepsilon) \circ \Delta \,\,=\,\, \coring{C} \,\,=\,\,
(\varepsilon\tensor{A}\coring{C}) \circ \Delta.$$ A \emph{right
$\coring{C}$-comodule} is a pair $(M,\varrho^M)$ with $M$ a right
$A$-module and $\varrho^M: M \to M\tensor{A}\coring{C}$ a right
$A$-linear map (called right $\coring{C}$-coaction) satisfying two
equalities $(\varrho^M \tensor{A} \coring{C}) \circ \varrho^M
\,=\, (M\otimes_A \Delta) \circ \varrho^M$, and $(M\otimes_A
\varepsilon) \circ \varrho^M \,=\,M$. A \emph{morphism}  of right
$\coring{C}$-comodules $f:(M,\varrho^M) \to (M',\varrho^{M'})$ is
a right $A$-linear map $f:M\to M'$ which is compatible with
coactions: $f \circ \varrho^{M'} \,=\, (f\tensor{A} \coring{C})
\circ \varrho^{M}$ ($f$ is right $\coring{C}$-colinear). The
category of all right $\coring{C}$-comodules is denoted by
$\rcomod{\coring{C}}$. Left $\coring{C}$-comodules are
symmetrically defined, we use the Greek letter $\lambda^{-}$ to
denote their coactions. The category of (right) $\cc$-comodules
is not in general an abelian category, it has cokernels and
arbitrary direct sums which can be already computed in the
category of $A$-modules. However, if ${}_A\cc$ is a flat module, then
$\rcomod{\cc}$ becomes a Grothendieck category (see
\cite{Kaoutit/Gomez/Lobillo:2004c}). Let $\dd$ be a $B$-coring, a
$(\cc,\dd)$-\emph{bicomodule} is a three-tuple
$(M,\varrho^M,\lambda^M)$ consisting of an $(A,B)$-bimodule ($M
\in \Bimod{A}{B}$) and $A-B$-bilinear maps $\varrho^M: M \to M
\otimes_B\coring{D}$ and $\lambda^M: M \to \coring{C}\otimes_A M$
such that $(M, \varrho^M)$ is right $\coring{D}$-comodule and
$(M,\lambda^M)$ is left $\coring{C}$-comodule with compatibility
condition $(\coring{C}\otimes_A \varrho^M) \circ \lambda^M  \,=\,
(\lambda^M \otimes_B \coring{D}) \circ \varrho^M$. A
\emph{morphism of bicomodules} is a left and right colinear map.
The category of all $(\coring{C},\coring{D})$-bicomodule is
denoted by $\Bicomod{\coring{C}}{\coring{D}}$. Obviously any ring
$A$ can be endowed with a trivial structure of an $A$-coring with
comultiplication the isomorphism $A\cong A\tensor{A}A$ (i.e.
$\iota_A$) and counit the identity $A$. In this way an
$(A,\dd)$-bicomodule is just a right $\dd$-comodule
$(M,\varrho^M)$ whose underlying module $M$ is an $A-B$-bimodule
and whose coaction $\varrho^M$ is an $A-B$-bilinear map.

Recall from \cite{Brzezinski/Kaoutit/Gomez:2006} (see
\cite{Lack/Street:2002} for general notions), that the bicategory
$\Sf{REM(Bim)}$ is defined by the following data
\begin{enumerate}[$\bullet$]
\item $0$-\emph{cells}. Are all corings $(\coring{C}:A)$ (i.e.,
$\coring{C}$ is an $A$-coring)

\item $1$-\emph{cells}. From $(\coring{D}:B)$ to $(\coring{C}:A)$
are pairs $(M,\mmm)$ where $M \in \Bimod{A}{B}$ and $\mmm:
\coring{C}\tensor{A}M \to M\tensor{B}\coring{D}$ is an
$A-B$-bilinear map satisfying
\begin{equation}\label{1-cell}
\begin{array}{c}
(M\tensor{B}\varepsilon_{\dd}) \circ \mmm \,\, =\,\,
\varepsilon_{\cc} \tensor{A}M,  \\ \\ (\mmm\tensor{B}\dd) \circ
(\coring{C}\tensor{A}\mmm) \circ (\Delta_{\coring{C}}\tensor{A}M)
\,\,=\,\, (M\tensor{B}\Delta_{\dd}) \circ \mmm
\end{array}
\end{equation}
where the first equality was token up to the isomorphism
$\iota_M$. The identity $1$-cell of $(\cc:A)$ is given by the pair
$\II_{(\cc:A)}\,=\,(A,\iota_{\cc}^{-1})$.

\item $2$-\emph{cells}. (in their reduced form) $\alpha: (M,\mmm) \longrightarrow
(M',\mmm')$ is an $A-B$-bilinear map $\alpha: \cc\tensor{A}M \to
M'$ satisfying
\begin{equation}\label{2-cell}
\mmm' \circ (\cc\tensor{A}\alpha) \circ (\Delta_{\cc}\tensor{A}M)
\,\,=\,\, (\alpha \tensor{B}\dd) \circ (\cc\tensor{A}\mmm) \circ
(\Delta_{\cc} \tensor{A}M).
\end{equation}
\end{enumerate}

Let $(M,\mmm)$ be a $1$-cell from $(\dd:B)$ to $(\cc:A)$, and
$(W,\www)$ a $1$-cell from $(\ee:C)$ to $(\dd:B)$. The horizontal
multiplication is defined by  $$(M,\mmm) (W,\www)\,\,=\,\,
\lr{M\tensor{B}W, (M\tensor{B}\www) \circ (\mmm\tensor{B}W)}$$ If
$\alpha: \cc\tensor{A}M \to M'$ and $\beta:\dd\tensor{B}W \to W'$
are two $2$-cells, then the vertical multiplication
$\alpha\beta:\cc\tensor{A}M\tensor{B}W \to M'\tensor{B}W'$, is
given by the rule
$$\alpha \beta\,\,=\,\, (M'\tensor{B}\beta) \circ
(\mmm'\tensor{B}W) \circ (\cc\tensor{A}\alpha\tensor{B}W) \circ
(\Delta_\cc \tensor{A}M\tensor{B}W)$$ The composition of $2$-cells
$\alpha: \cc\tensor{A}M \to M'$ and $\alpha':\cc\tensor{A}M' \to
M''$ is given by
$$\alpha' \boldsymbol{\circ} \alpha\,=\, \alpha'
\circ (\cc\tensor{A}\alpha) \circ (\Delta_\cc \tensor{A}M).$$
The left and right identity functors are given as follows. Consider
$(M,\mmm)$ any $1$-cell  from $(\dd:B)$ to $(\cc:A)$, the left
identity multiplication $\lll_{(M,\mmm)}: \II_{(\cc:\,A)}(M,\mmm)
\to (M,\mmm)$ is defined as the composition
$$\xymatrix{\lll_{(M,\mmm)}:\,\,\cc\tensor{A}A\tensor{A}M
\ar@{->}^-{\varepsilon\tensor{}M}[r] & A\tensor{A}A\tensor{A}M
\ar@{->}^-{{\bf l}}[r] & A\tensor{A}M \ar@{>}^-{{\bf l}}[r] & M
}$$ while the right identity multiplication $\rrr_{(M,\mmm)}:
(M,\mmm) \II_{(\dd:\,B)} \to (M,\mmm)$ is defined as the
composition
$$\xymatrix@C=40pt{\rrr_{(M,\mmm)}:\,\,\cc\tensor{A}M\tensor{B}B
\ar@{->}^-{\varepsilon\tensor{}M\tensor{}B}[r] &
A\tensor{A}M\tensor{B}B \ar@{->}^-{{\bf l}}[r] & M\tensor{B}B
\ar@{>}^-{{\bf r}}[r] & M. }$$

Let $(M,\varrho^M)$ be an $(A,\dd)$-bicomodule, and consider the
following composed map
$$\xymatrix@C=50pt{\mmm: \cc\tensor{A}M
\ar@{->}^-{\cc\tensor{}\varrho^M}[r] & \cc\tensor{A}M\tensor{B}\dd
\ar@{->}^-{\varepsilon_\cc\tensor{}M\tensor{}\dd}[r] &
A\tensor{A}M\tensor{B}\dd \,\,\cong\,\, M\tensor{B}\dd}$$

\begin{lemma}\label{bicomd-cell}
The pair $(M,\mmm)$ is a $1$-cell from $(\dd:B)$ to $(\cc:A)$ in
the bicategory $\Sf{REM(Bim)}$.
\end{lemma}
\begin{proof}
By definition $\mmm$ is an $A-B$-bilinear map. Let us show that
$\mmm$ satisfies the equations of \eqref{1-cell}. First, we have
\begin{eqnarray*}
  (M \tensor{B}\varepsilon_\dd) \circ \mmm &=& (M \tensor{B}\varepsilon_\dd)
  \circ (\varepsilon_\cc \tensor{A}M\tensor{B}\dd) \circ (\cc\tensor{A}\varrho^M) \\
   &=& (\varepsilon_\cc \tensor{A}M) \circ (\cc\tensor{A}M\tensor{B}\varepsilon_\dd)
    \circ (\cc\tensor{A} \varrho^M)  \\
   &=& (\varepsilon_\cc \tensor{A}M)
\end{eqnarray*}
and secondly, on one hand we have
\begin{eqnarray*}
  (M\tensor{B}\Delta_\dd) \circ \mmm &=& (M\tensor{B}\Delta_\dd) \circ
  (\varepsilon_\cc \tensor{A}M\tensor{B}\cc) \circ (\cc\tensor{A}\varrho^M) \\
   &=& (\varepsilon_\cc\tensor{A}M\tensor{B}\dd\tensor{B}\dd) \circ (\cc\tensor{A}M\tensor{B}\Delta_\dd)
   \circ (\cc\tensor{A}\varrho^M) \\
   &=& (\varepsilon_\cc\tensor{A}M\tensor{B}\cc\tensor{B}\dd) \circ
   \lr{\cc\tensor{A}\lr{(M\tensor{B}\Delta_\dd) \circ \varrho^M }} \\
   &=& (\varepsilon_\cc\tensor{A}M\tensor{B}\cc\tensor{B}\dd) \circ
   \lr{\cc\tensor{A}\lr{(\varrho^M\tensor{B}\dd) \circ \varrho^M }} \\
   &=& (\varepsilon_\cc\tensor{A}M\tensor{B}\dd\tensor{B}\dd) \circ
   (\cc\tensor{A}\varrho^M\tensor{B}\dd) \circ (\cc\tensor{A}
   \varrho^M)\\ &=& (\mmm\tensor{B}\dd) \circ (\cc\tensor{A}\varrho^M)
\end{eqnarray*}
and on the other hand we have
\begin{eqnarray*}
 (\mmm\tensor{B}\dd) \circ (\cc\tensor{A}\mmm) \circ (\Delta_\cc \tensor{A}M)  &=&
  \lr{\lr{(\varepsilon_\cc\tensor{A}M\tensor{B}\dd) \circ (\cc\tensor{A}\varrho^M)}\tensor{B}\dd}
  \\ &\,\,&\; \circ
  \lr{ \cc\tensor{A}\lr{(\varepsilon_\cc\tensor{A}M\tensor{B}\dd) \circ (\cc\tensor{A}\varrho^M)}}
  \circ (\Delta_\cc \tensor{A}M)\\
   &=& (\varepsilon_\cc\tensor{A}M\tensor{B}\dd\tensor{B}\dd) \circ (\cc\tensor{A}\varrho^M \tensor{B}\dd)\circ
  (\cc\tensor{A}\varepsilon_\cc\tensor{A}M\tensor{B}\dd) \\ &\,\,&\; \circ (\cc\tensor{A}\cc\tensor{A}\varrho^M)
   \circ (\Delta_\cc \tensor{A}M) \\
   &=& (\varepsilon_\cc\tensor{A}M\tensor{B}\dd\tensor{B}\dd) \circ (\cc\tensor{A}\varrho^M \tensor{B}\dd) \circ
  (\cc\tensor{A}\varepsilon_\cc\tensor{A}M\tensor{B}\dd) \\ &\,\,&\;
  \circ (\Delta_\cc\tensor{A}M\tensor{B}\dd) \circ (\cc\tensor{A}\varrho^M) \\
   &=& (\varepsilon_\cc\tensor{A}M\tensor{B}\dd\tensor{B}\dd) \circ (\cc\tensor{A}\varrho^M \tensor{B}\dd) \circ
      (\cc\tensor{A}\varrho^M) \\
   &=& (\mmm\tensor{B}\dd) \circ (\cc\tensor{A}\varrho^M)
\end{eqnarray*}
Therefore, $(M\tensor{B}\Delta_\dd) \circ \mmm\,=\,
(\mmm\tensor{B}\dd) \circ (\cc\tensor{A}\mmm) \circ (\Delta_\cc
\tensor{A}M)$.
\end{proof}

The bicategory $\Sf{W(REM(Bim))}$ is then defined by the following
data:

\begin{enumerate}[$\bullet$]
\item $0$-\emph{cells}. Are all corings $(\cc:A)$.

\item $1$-\emph{cells}. From $(\coring{D}:B)$ to $(\coring{C}:A)$,
are four-tuples $\Gamma=\lr{(M,\mmm),(N,\nnn), \eta ,\rho}$ with
$M \in \Bimod{A}{B}$, $N \in \Bimod{B}{A}$, bilinear maps $\mmm:
\cc\tensor{A}M \to M\tensor{B}\dd$, $\nnn: \dd\tensor{B}N \to
N\tensor{A}\cc$ satisfying \eqref{1-cell}, and two $2$-cells,
$\eta: \cc\tensor{A}M\tensor{B}N \to A$, and
$\rho:\dd\tensor{B}N\tensor{A}M \to B$ satisfying both them the
following equalities (the first two ones are up to the natural
isomorphism $\iota_{-}$)
\begin{multline}\label{I-1}
(\eta\tensor{A}\cc) \circ (\cc\tensor{A}M\tensor{B}\nnn)
  \circ (\cc\tensor{A}\mmm\tensor{B}N) \circ (\Delta_\cc
  \tensor{A}M\tensor{B}N)\,=\, (\cc\tensor{A}\eta) \circ (\Delta_\cc \tensor{A}M\tensor{B}N),
\end{multline}
\begin{multline}\label{I-2}
(\rho \tensor{B} \dd) \circ (\dd\tensor{B}N\tensor{A}\mmm) \circ
(\dd\tensor{B}\nnn\tensor{A}M) \circ
(\Delta_\dd\tensor{B}N\tensor{A}M) \,=\, (\dd\tensor{B}\rho)
\circ (\Delta_\dd\tensor{B}N\tensor{A}M),
\end{multline}
\begin{eqnarray}
  \iota_M \circ (\eta \tensor{A}M)  &=& (M\tensor{B}\rho) \circ (\mmm\tensor{B}N\tensor{A}M), \label{II-1}\\
   \iota_N \circ (\rho\tensor{B}N) &=& (N\tensor{A}\eta) \circ
   (\nnn\tensor{A}M\tensor{B}N). \label{II-2}
\end{eqnarray}
The identity $1$-cell of $(\cc:A)$ is given by the four-tuple
$\lr{(A,\iota_\cc^{-1}),\,(A,\iota_\cc^{-1}),\,\rrr_{(\cc:\,
A)},\,\lll_{(\cc:\,A)}}$

\item $2$-\emph{cells}. Are pairs
$(\alpha,\beta):\Gamma=\lr{(M,\mmm),(N,\nnn), \eta ,\rho} \to
\Gamma'=\lr{(M',\mmm'),(N',\nnn'), \eta' ,\rho'}$ where $\alpha:
\cc\tensor{A}M \to M'$ and $\beta:\dd\tensor{B}N \to N'$ are
$2$-cells satisfying the following equalities
\begin{eqnarray}
  \eta &=& \eta' \circ (\cc\tensor{A}M'\tensor{B}\beta) \circ (\cc\tensor{A}\mmm'\tensor{B}N) \circ
  (\cc\tensor{A}\cc\tensor{A}\alpha\tensor{B}N)
  \\ &\,\,& \; \circ\,\, (\cc\tensor{A}\Delta_\cc\tensor{A}M\tensor{B}N)
  \circ ( \Delta_\cc\tensor{A}M\tensor{B}N) \nonumber \\
  \rho &=& \rho' \circ (\dd\tensor{B}N'\tensor{A}\alpha)
  \circ (\dd\tensor{B}\nnn'\tensor{A}M) \circ
  (\dd\tensor{B}\dd\tensor{B}\beta\tensor{A}M) \\ &\,\,&\; \circ
  \,\, (\dd\tensor{B}\Delta_\dd\tensor{B}N\tensor{A}M) \circ
  (\Delta_\dd\tensor{B}N\tensor{A}M)\nonumber
\end{eqnarray}
\end{enumerate}

\begin{rem}
Recall that a (classical) Morita context between two rings $A$ and $B$
(or from $B$ to $A$) consists of two bimodules $P \in \Bimod{A}{B}$
and $Q \in \Bimod{B}{A}$ together with bilinear maps
$\varphi:P\tensor{B}Q \to A$ and $\psi: Q\tensor{A}P \to B$
satisfying $$\,\, y\varphi(x,y')\,\,=\,\, \psi(y,x)y', \qquad
x\psi(y, x')\,\,=\,\,\varphi(x,y)x'$$ for every $x,x' \in P$ and $y,
y' \in Q$. This  of course says exactly that the four-tuple
$(P,Q,\varphi,\psi)$ is a WRMC from $B$ to $A$ in the bicategory
$\Sf{W(Bim)}$. This can also be recognized in the bicategory $\Sf{W(REM(Bim))}$ as
follows. As was mentioned before we can consider trivially any ring
$A$ as an $A$-coring $(A:A)$. That is we can consider $(A:A)$ as a
$0$-cell in $\Sf{REM(Bim)}$ for every ring $A$. So let $(A:A)$ and
$(B:B)$ two trivial $0$-cells in $\Sf{REM(Bim)}$. The category of
all $1$-cells from $(B:B)$ to $(A:A)$ in this bicategory is clearly,
using the natural isomorphism $\iota_{-}$, isomorphic to the
category of bimodules $\Bimod{A}{B}$. Now, consider
$\Gamma=((M,\iota_M), (N,\iota_N), \eta,\rho)$ any $1$-cell from
$(B:B)$ to $(A:A)$ in the bicategory $\Sf{W(REM(Bim))}$. This means
that $M \in \Bimod{A}{B}$, $N \in \Bimod{B}{A}$, and
$\eta:A\tensor{A}M\tensor{B}N \to A$, $\rho:B\tensor{B}N\tensor{A}M
\to B$ are bilinear maps satisfying equations
\eqref{I-1}-\eqref{II-2}. Define in an obvious way
$\widehat{\eta}:M\tensor{B}N \to A$ and $\widehat{\rho}: N\tensor{A}M
\to B$. It is easily seen that $(M,N,\widehat{\eta},
\widehat{\rho})$ is a Morita context between $A$ and $B$. Conversely
any WRMC from $(B:B)$ to $(A:A)$ in $\Sf{W(REM(Bim))}$ is deduced
from a classical Morita context from $B$ to $A$.
\end{rem}

Fix an arbitrary $1$-cell $\Gamma=\lr{(M,\mmm), (N,\nnn), \eta,
\rho}$ in ${}_{(\cc:\,A)}\Sf{W(REM(Bim))}_{(\dd:\,B)}$. Since
$(M,\mmm)$ and $(N,\nnn)$ are a $1$-cells in $\Sf{REM(Bim)}$, we can
associate by \cite[4.2]{Brzezinski/Kaoutit/Gomez:2006}, the (right)
push-out functors $\M:\rcomod{\cc} \to \rcomod{\dd}$ and $ \Nn :
\rcomod{\dd} \to \rcomod{\cc}$, explicitly acting on objects by
\begin{eqnarray}\label{push-out}
  \M(X,\varrho^X) &=& \lr{X\tensor{A}M,\, \varrho^{X\tensor{A}M}=
(X\tensor{A}\mmm) \circ (\varrho^X\tensor{A}M)} \\
 \Nn(Y,\varrho^Y)  &=& \lr{Y\tensor{B}N,\, \varrho^{Y\tensor{B}N}=
(Y\tensor{B}\nnn) \circ (\varrho^Y\tensor{B}N)}, \nonumber
\end{eqnarray}
for every right $\cc$-comodule $(X,\varrho^X)$ and every right
$\dd$-comodule $(Y,\varrho^Y)$. The actions on morphisms are
obvious.

On the other hand, we can define the following right linear maps
$$ \xymatrix@R=8pt@C=50pt{\widetilde{\eta}_X: \, X\tensor{A}M\tensor{B}N
\ar@{->}^-{\varrho^X\tensor{}M\tensor{}N}[r] &
X\tensor{A}\cc\tensor{A}M\tensor{B}N \ar@{->}^-{X\tensor{}\eta}[r]
& X\tensor{A}A \cong X \\ \widetilde{\rho}_Y: \,
Y\tensor{B}N\tensor{A}M
\ar@{->}^-{\varrho^Y\tensor{}N\tensor{}M}[r] &
Y\tensor{B}\dd\tensor{B}N\tensor{A}M
\ar@{->}^-{Y\tensor{}\varrho}[r] & Y\tensor{B}B \cong Y}$$ for
every pair of right $\cc$-comodule $(X,\varrho^X)$ and right
$\dd$-comodule $(Y,\varrho^Y)$. We consider $\Nn
\M(X)\,=\,X\tensor{A}M\tensor{B}N$ and $\M
\Nn(Y)\,=\,Y\tensor{B}N\tensor{A}M$, respectively, as right
$\cc$-comodule and as right $\dd$-comodule using coactions defined
in \eqref{push-out}.

\begin{prop}\label{Morita-comod}
Keeping the above notations, we have
\begin{enumerate}[(i)]
\item The right linear maps $\widetilde{\eta}_X$ and
$\widetilde{\rho}_Y$ are, respectively, right $\cc$-colinear and
right $\dd$-colinear.

\item They define a natural transformation $\widetilde{\eta}_{-}:
\Nn\M \to \id_{\rcomod{\cc}}$ and $\widetilde{\rho}: \M \Nn \to
\id_{\rcomod{\dd}}$ satisfying the following conditions
\begin{eqnarray*}
  \widetilde{\eta}_{\Nn(Y)} &=& \Nn \widetilde{\rho}_{Y},\quad \text{for\,\,every}\,\, (Y,\varrho^Y) \in \rcomod{\dd},  \\
  \widetilde{\rho}_{\M(X)} &=& \M \widetilde{\eta}_{X},\quad \text{for\,\,every}\,\, (X ,\varrho^X)\in
  \rcomod{\cc}.
\end{eqnarray*}
That is $(\M,\Nn,\widetilde{\eta},\widetilde{\varrho})$ is right
wide Morita context, in the sense of \cite{Castano/Gomez:1995},
between the categories $\rcomod{\cc}$ and $\rcomod{\dd}$ (i.e. a WRMC in the $2$-category $\Sf{Cat}$ with right exact functors).
\end{enumerate}
\end{prop}
\begin{proof}
$(i)$. We only prove the colinearity of $\widetilde{\eta}_X$, a
similar arguments are used for $\widetilde{\rho}_Y$. By definition
$\widetilde{\eta}_X$ is right $A$-linear, and we have
\begin{eqnarray*}
  \varrho^X \circ \widetilde{\eta}_X &=& \varrho^X \circ (X\tensor{A}\eta) \circ (\varrho^X \tensor{A}M\tensor{B}N) \\
   &=& (X\tensor{A}\cc\tensor{A}\eta) \circ
   (\varrho^X\tensor{A}\cc\tensor{A}M\tensor{B}N)\circ
   (\varrho^X\tensor{A}M\tensor{B}N)\\
   &=& (X\tensor{A}\cc\tensor{A}\eta) \circ
   \lr{\lr{(\varrho^X\tensor{A}\cc) \circ \varrho^X}\tensor{A}M\tensor{B}N} \\
   &=& (X\tensor{A}\cc\tensor{A}\eta) \circ
   \lr{\lr{(X\tensor{A}\Delta_\cc) \circ \varrho^X}\tensor{A}M\tensor{B}N} \\
   &=& (X\tensor{A}\cc\tensor{A}\eta) \circ
   (X\tensor{A}\Delta_\cc\tensor{A}M\tensor{B}N) \circ (\varrho^X\tensor{A}M\tensor{B}N) \\
   &=& \lr{X\tensor{A}\lr{ (\cc\tensor{A}\eta) \circ (\Delta_\cc\tensor{A}M\tensor{B}N)}}
   \circ (\varrho^X\tensor{A}M\tensor{B}N) \\
   &=& \lr{X\tensor{A}\lr{ (\eta\tensor{A}\cc) \circ
   (\cc\tensor{A}M\tensor{B}\nnn) \circ
   (\cc\tensor{A}\mmm\tensor{B}N) \circ
   (\Delta_\cc\tensor{A}M\tensor{B}N)}} \\ &\,\,&\,\,\,
   \circ\, (\varrho^X\tensor{A}M\tensor{B}N),\text{\,\,by}\,\, \eqref{I-1} \\
   &=& (X\tensor{A}\eta\tensor{A}\cc) \circ
   (X\tensor{A}\cc\tensor{A}M\tensor{B}\nnn) \circ
   (X\tensor{A}\cc\tensor{A}\mmm\tensor{B}N) \\ &\,\,&\,\,\, \circ\,
   (X\tensor{A}\Delta_\cc\tensor{A}M\tensor{B}N)
   \circ (\varrho^X\tensor{A}M\tensor{B}N) \\
   &=& (X\tensor{A}\eta\tensor{A}\cc) \circ
   (X\tensor{A}\cc\tensor{A}M\tensor{B}\nnn) \circ
   (X\tensor{A}\cc\tensor{A}\mmm\tensor{B}N) \\ &\,\,&\,\,\, \circ\,
   (\varrho^X\tensor{A}\cc\tensor{A}M\tensor{B}N)
   \circ (\varrho^X\tensor{A}M\tensor{B}N) \\
   &=& (X\tensor{A}\eta\tensor{A}\cc) \circ
   (X\tensor{A}\cc\tensor{A}M\tensor{B}\nnn) \circ
   (\varrho^X\tensor{A}M\tensor{B}\dd\tensor{B}N) \\ &\,\,&\,\,\, \circ\,
   (X\tensor{A}\mmm\tensor{A}M\tensor{B}N)
   \circ (\varrho^X\tensor{A}M\tensor{B}N) \\
   &=& (X\tensor{A}\eta\tensor{A}\cc) \circ
   (\varrho^X\tensor{A}M\tensor{B}N\tensor{A}\cc) \circ
   (X\tensor{A}M\tensor{B}\nnn) \\ &\,\,&\,\,\, \circ\,
   (X\tensor{A}\mmm\tensor{A}M\tensor{B}N)
   \circ (\varrho^X\tensor{A}M\tensor{B}N) \\
   &=& (\widetilde{\eta}_X\tensor{A}\cc) \circ
   \varrho^{X\tensor{A}M\tensor{B}N},
\end{eqnarray*}
where $\varrho^{X\tensor{B}M\tensor{A}N}\,=\,
(X\tensor{B}M\tensor{A}\nnn) \circ
(X\tensor{B}\mmm\tensor{B}M\tensor{A}N)\circ
(\varrho^X\tensor{B}M\tensor{A}N)$ is the right coaction of
$\Nn\M(X,\varrho^X)$ as defined in \eqref{push-out}.\\ $(ii)$. Let
$f:(X,\varrho^X) \to (X',\varrho^{X'})$ be a right $\cc$-colinear
map, then
\begin{eqnarray*}
  f \circ \widetilde{\eta}_X &=& f \circ (X\tensor{A}\eta) \circ (\varrho^X \tensor{A}M\tensor{B}N) \\
   &=& (X'\tensor{A}\eta) \circ (f\tensor{A}\cc\tensor{A}M\tensor{B}N)
   \circ (\varrho^X \tensor{A}M\tensor{B}N), \,\,f\text{\,\, is \,\, linear} \\
   &=& (X'\tensor{A}\eta) \circ \lr{\lr{(f\tensor{A}\cc) \circ
   \varrho^X}\tensor{A}M\tensor{B}N} \\
   &=& (X'\tensor{A}\eta) \circ (\varrho^{X'}\tensor{A}M\tensor{B}N)
   \circ (f \tensor{A}M\tensor{B}N),\,\,f\text{\,\, is \,\, colinear}  \\
   &=& \widetilde{\eta}_{X'} \circ \Nn\M(f).
\end{eqnarray*} Hence $\widetilde{\eta}_{-}$ is natural. Similarly we check that
$\widetilde{\rho}_{-}$ is also natural. Given now any right
$\dd$-comodule $(Y,\varrho^Y)$, we know that $\varrho^{\Nn(Y)}=
(N\tensor{A}\nnn) \circ (\varrho^Y\tensor{B}N)$, so
\begin{eqnarray*}
  \widetilde{\eta}_{\Nn(Y)} &=& (Y\tensor{B}N\tensor{A}\eta) \circ (Y\tensor{B}\nnn\tensor{A}M\tensor{B}N)
  \circ (\varrho^Y\tensor{B}N\tensor{A}M\tensor{B}N) \\
   &=& \lr{Y\tensor{B}\lr{(N\tensor{A}\eta) \circ (\nnn\tensor{A}M\tensor{B}N)}}
   \circ (\varrho^Y\tensor{B}N\tensor{A}M\tensor{B}N) \\
   &=& (Y\tensor{B}\rho\tensor{B}N) \circ (\varrho^Y\tensor{B}N\tensor{A}M\tensor{B}N),
   \,\, \text{by\,\,} \eqref{II-2}\,\,(\text{up \,\,to}\,\,\iota_{-})  \\
   &=& \Nn\widetilde{\rho}_Y,
\end{eqnarray*} which gives the first condition of the stated item $(ii)$. The
second one is similarly proved. Since the push-out functors preserve
cokernels, we have thus showed that
$(\M,\Nn,\widetilde{\eta},\widetilde{\rho})$ is right Morita context
between $\rcomod{\dd}$ and $\rcomod{\cc}$ in the sense of
\cite{Castano/Gomez:1995}.
\end{proof}

\begin{rem}
A more conceptual proof of Proposition \ref{Morita-comod} can be given using Proposition \ref{Morph-bicatg}, in the following way. One can show that the push-out defines a homomorphism of bicategories $\Scr{P}: \Sf{REM(Bim)} \to \Sf{Cat}$, sending any $0$-cell $(\cc:A)$ to its category of right $\cc$-comodules $\rcomod{\cc}$. Therefore, there is a homomorphism of bicategories $\Sf{W}(\Scr{P}):\Sf{W(REM(Bim))} \to \Sf{W}(\Sf{Cat})$. Since the push-out functor preserves cokernels and direct sums, we immediately deduce the claim of Proposition \ref{Morita-comod}.\\
On the other hand, push-out functors are in fact a lifting functors. Precisely, given a bimodule ${}_AM_B$, we know (see \cite[p. 256]{Lack/Street:2002} and the dual version of \cite[Lemma 1]{Johnstone:1975}) that there is a 1-1 correspondence between
\begin{enumerate}[(i)]
\item $A-B$-bilinear maps $\mmm: \cc\tensor{A}M \to M\tensor{B}\dd$ satisfying equation \eqref{1-cell};
\item Lifting functors $\mathcal{M}: \rcomod{\cc} \to \rcomod{\dd}$ of the functor $-\tensor{A}M : \rmod{A} \to \rmod{B}$, that is, functors $\mathcal{M}$ rendering commutative the following diagram
    $$
    \xymatrix{ \rcomod{\cc} \ar@{-->}^-{\mathcal{M}}[rr] \ar@{->}_-{\mathcal{O}_{\cc}}[d] & & \rcomod{\dd} \ar@{->}^-{\mathcal{O}_{\dd}}[d] \\ \rmod{A} \ar@{->}_-{-\tensor{A}M}[rr] & &  \rmod{B} }
    $$ where $\Oo_{\cc}$ and $\Oo_{\dd}$ are the forgetful functors.
\end{enumerate}
In this way, it is posible to analyse the converse of Proposition \ref{Morita-comod}. The forthcoming proposition  discusses a special case of this analysis.
\end{rem}

Suppose that a wide right Morita context in the sense of
\cite{Castano/Gomez:1995} between $\rcomod{\dd}$ and $\rcomod{\cc}$
is given. This consists of pair of covariant additive functors
$\xymatrix{ \Mm:\,\rcomod{\dd} \ar@<0,5ex>[r] & \rcomod{\cc}:\,\Nn
\ar@<0,5ex>[l]}$ which preserve cokernels, and two natural
transformations $\widetilde{\eta}: \Mm\Nn \to \id_{\rcomod{\cc}}$
and $\widetilde{\rho}: \Nn\Mm \to \id_{\rcomod{\dd}}$ with
compatible conditions as stated in Lemma \ref{Morita-comod}(ii). Assume that $\Mm$ and $\Nn$ preserve direct sums, and they are lifting functors for some $(A,\dd)$-bicomodule $(M,\varrho^M)$
and some $(B,\cc)$-bicomodule $(N,\varrho^N)$, respectively (here
$A$ and $B$ are considered trivially as an $A$-coring and
$B$-coring). In this case, we are assuming that there are commutative
diagrams
\begin{equation}\label{Hyp-I}
\xymatrix{\rcomod{\cc} \ar@{->}^-{\Mm}[rr] \ar@{->}_-{\Oo_\cc}[dr] &
& \rcomod{\dd} \\ & \rmod{A} \ar@{-->}_-{-\tensor{A}M}[ur] & }
\qquad \xymatrix{\rcomod{\dd} \ar@{->}^-{\Nn}[rr]
\ar@{->}_-{\Oo_\dd}[dr] & & \rcomod{\cc} \\ & \rmod{B}
\ar@{-->}_-{-\tensor{B}N}[ur] & }
\end{equation}
We know from Lemma \ref{bicomd-cell}, that $(M ,\mmm)$ is a $1$-cell in
$\Sf{REM(Bim)}$ from $(\dd:B)$ to $(\cc:A)$, where
$$\xymatrix@C=40pt{\mmm:\,\cc\tensor{A}M
\ar@{->}^-{\cc\tensor{}\varrho^M}[r] & \cc\tensor{A}M\tensor{B}\dd
\ar@{->}^-{\varepsilon_\cc\tensor{}M\tensor{}\dd}[r] &
A\tensor{A}M\tensor{B}\dd \cong M\tensor{B}\dd.}$$ Similarly
$(N,\nnn)$ is a $1$-cell from $(\cc:A)$ to $(\dd:B)$ in
$\Sf{REM(Bim)}$ with $$\xymatrix@C=40pt{\nnn:\,\dd\tensor{B}N
\ar@{->}^-{\dd\tensor{}\varrho^N}[r] & \dd\tensor{B}N\tensor{A}\cc
\ar@{->}^-{\varepsilon_\dd\tensor{}N\tensor{}\cc}[r] &
B\tensor{B}N\tensor{A}\cc \cong N\tensor{A}\cc}$$

\begin{prop}\label{Reciproco}
Let $\Mm$, $\Nn$ be the functors defined in \eqref{Hyp-I}, and
$(M,\mmm)$, $(N,\nnn)$, $\widetilde{\eta}$ and $\widetilde{\rho}$
as above. Then $\Gamma\,=\,\lr{ (M,\mmm),(N,\nnn), \eta,\rho }$ is
a $1$-cell in ${}_{(\cc:\,A)}\Sf{W(REM(Bim))}_{(\dd:\,B)}$, where
$\eta\,=\, \varepsilon_\dd \circ \widetilde{\eta}_\dd$ and
$\rho\,=\, \varepsilon_\cc \circ \widetilde{\rho}_\cc$.
\end{prop}
\begin{proof}
We need to show that the stated maps $\eta$ and $\rho$ satisfy
equations \eqref{I-1}-\eqref{II-2}.

Since $\widetilde{\eta}_{-}$ and $\widetilde{\rho}_{-}$ are natural
transformations, it is clearly seen that, for every ring $R$ and
every pair of $(R,\dd)$-bicomodule $(Y,\varrho^Y)$ and
$(R,\cc)$--bicomodule $(X,\varrho^X)$, we have
$\widetilde{\eta}_{Y}: Y\tensor{B}M\tensor{A}N \to Y$ and
$\widetilde{\rho}_X: X\tensor{A}N\tensor{B}M \to X$ are,
respectively, $R-B$-bilinear and $R-A$-bilinear maps. Furthermore,
we can prove, using free presentations of right modules, that
\begin{eqnarray*}
  \widetilde{\eta}_{Y\tensor{B}\dd} &=& Y\tensor{B}\widetilde{\eta}_{\dd},
  \quad\text{for\,\, every \,\, right \,\,module}\,\, Y \in \rmod{B} \\
  \widetilde{\rho}_{X\tensor{A}\cc} &=&
  X\tensor{A}\widetilde{\rho}_\cc,\quad\text{for\,\, every \,\,right\,\, module}\,\, X \in
  \rmod{A}.
\end{eqnarray*}
Keeping this in mind, we compute
\begin{eqnarray*}
  & & (\eta\tensor{A}\cc)   \circ (\cc\tensor{A}M\tensor{B}\nnn) \circ (\cc\tensor{A}\mmm\tensor{B}N)
  \circ (\Delta_\cc\tensor{A}M\tensor{B}\dd\tensor{B}N)  \\
   &=& (\varepsilon_\cc\tensor{A}\cc) \circ (\widetilde{\eta}_\cc\tensor{A}\cc)
   \circ (\cc\tensor{A}M\tensor{B}\nnn) \circ (\cc\tensor{A}\mmm\tensor{B}N)
   \\ &\,\,& \circ\, (\Delta_\cc\tensor{A}M\tensor{B}\dd\tensor{B}N)  \\
   &=& (\varepsilon_\cc\tensor{A}\cc) \circ (\widetilde{\eta}_\cc\tensor{A}\cc)
   \circ (\cc\tensor{A}M\tensor{B}\varepsilon_\dd\tensor{B}N\tensor{A}\cc)
   \circ (\cc\tensor{A}M\tensor{B}\dd\tensor{B}\varrho^N) \\
   &\,\,& \circ \, (\cc\tensor{A}\varepsilon_\cc\tensor{A}M\tensor{B}\dd\tensor{B}N)
   \circ (\cc\tensor{A}\cc\tensor{A}\varrho^M\tensor{B}N) \circ (\Delta_\cc\tensor{A}M\tensor{B}\dd\tensor{B}N)  \\
   &=& (\varepsilon_\cc\tensor{A}\cc) \circ (\widetilde{\eta}_\cc\tensor{A}\cc)
   \circ (\cc\tensor{A}M\tensor{B}\varepsilon_\dd\tensor{B}N\tensor{A}\cc)
   \circ (\cc\tensor{A}M\tensor{B}\dd\tensor{B}\varrho^N) \\
   &\,\,& \circ \, (\cc\tensor{A}\varepsilon_\cc\tensor{A}M\tensor{B}\dd\tensor{B}N)
   \circ (\Delta_\cc \tensor{A}M\tensor{B}\dd\tensor{B}N) \circ (\cc\tensor{A}\varrho^M\tensor{B}N)   \\
   &=& (\varepsilon_\cc\tensor{A}\cc) \circ (\widetilde{\eta}_\cc\tensor{A}\cc)
   \circ (\cc\tensor{A}M\tensor{B}\varepsilon_\dd\tensor{B}N\tensor{A}\cc)
   \circ (\cc\tensor{A}M\tensor{B}\dd\tensor{B}\varrho^N) \\
   &\,\,& \circ \, (\cc\tensor{A}\varrho^M\tensor{B}N) \\
   &=& (\varepsilon_\cc\tensor{A}\cc) \circ (\widetilde{\eta}_\cc\tensor{A}\cc)
   \circ (\cc\tensor{A}M\tensor{B}\varepsilon_\dd\tensor{B}N\tensor{A}\cc)
   \circ (\cc\tensor{A}\varrho^M\tensor{B}N\tensor{A}\cc) \\
   &\,\,& \circ \, (\cc\tensor{A}M\tensor{B}\varrho^N) \\
   &=& (\varepsilon_\cc\tensor{A}\cc) \circ (\widetilde{\eta}_\cc\tensor{A}\cc)
   \circ (\cc\tensor{A}M\tensor{B}\varrho^N) \\
   &=& (\varepsilon_\cc\tensor{A}\cc) \circ \Delta_\cc \circ
   \widetilde{\eta}_\cc, \quad \widetilde{\eta}_\cc \text{\,\, is
   \,\, colinear} \\ &=& \widetilde{\eta}_\cc.
\end{eqnarray*}
On the other hand, we have
\begin{eqnarray*}
  (\cc\tensor{A}\eta) \circ (\Delta_\cc\tensor{A}\tensor{B}N)  &=&
  (\cc\tensor{A}\varepsilon_\cc) \circ (\cc\tensor{A}\widetilde{\eta}_\cc) \circ (\Delta_\cc\tensor{A}M\tensor{B}N) \\
   &=&  (\cc\tensor{A}\varepsilon_\cc) \circ
   \widetilde{\eta}_{\cc\tensor{A}\cc} \circ (\Delta_\cc\tensor{A}M\tensor{B}N) \\
   &=& (\cc\tensor{A}\varepsilon_\cc) \circ \Delta_\cc \circ
   \widetilde{\eta}_{\cc} \,\,=\,\, \widetilde{\eta}_{\cc}
\end{eqnarray*}
Therefore, $\eta$ satisfies \eqref{I-1}. Similarly, we prove that
$\rho$ satisfies \eqref{I-2}. Now, up to the isomorphism $\iota_M$,
we compute
\begin{eqnarray*}
  && (M\tensor{B}\varepsilon_\dd) \circ (M\tensor{B}\widetilde{\rho}_\dd) \circ
  (\varepsilon_\cc\tensor{A}M\tensor{B}\dd\tensor{B}N\tensor{A}M)
  \circ (\cc\tensor{A}\varrho^M\tensor{B}N\tensor{A}M) \\
  &=& (M\tensor{B}\varepsilon_\dd) \circ (\varepsilon_\cc\tensor{A}M\tensor{B}\dd) \circ
  (\cc\tensor{A}M\tensor{B}\widetilde{\rho}_\dd) \circ  (\cc\tensor{A}\varrho^M\tensor{B}N\tensor{A}M) \\
   &=& (M\tensor{B}\varepsilon_\dd) \circ (\varepsilon_\cc\tensor{A}M\tensor{B}\dd) \circ
  (\cc\tensor{A}\widetilde{\rho}_{M\tensor{B}\dd}) \circ  (\cc\tensor{A}\varrho^M\tensor{B}N\tensor{A}M) \\
   &=& (M\tensor{B}\varepsilon_\dd) \circ (\varepsilon_\cc\tensor{A}M\tensor{B}\dd) \circ
  \lr{\cc\tensor{A}\lr{ \widetilde{\rho}_{M\tensor{B}\dd} \circ  (\varrho^M\tensor{B}N\tensor{A}M)}}  \\
   &=& (M\tensor{B}\varepsilon_\dd) \circ (\varepsilon_\cc\tensor{A}M\tensor{B}\dd) \circ
  \lr{\cc\tensor{A}\lr{ \varrho^M \circ  \widetilde{\rho}_M}}  \\
   &=& (M\tensor{B}\varepsilon_\dd) \circ (\varepsilon_\cc\tensor{A}M\tensor{B}\dd) \circ
  (\cc\tensor{A}\varrho^M) \circ (\cc\tensor{A}\widetilde{\rho}_M) \\
   &=& (M\tensor{B}\varepsilon_\dd) \circ \varrho^M \circ (\varepsilon_\cc\tensor{A}M) \circ
  (\cc\tensor{A}\widetilde{\rho}_M) \\ &=& (\varepsilon_\cc\tensor{A}M) \circ
  (\cc\tensor{A}\widetilde{\rho}_M) \,\, =\,\, (\varepsilon_\cc\tensor{A}M) \circ
  \widetilde{\rho}_{\cc\tensor{A}M}\,\,=\,\, (\varepsilon_\cc\tensor{A}M) \circ
  (\widetilde{\eta}_\cc \tensor{A}M),
\end{eqnarray*}
which implies the equation \eqref{II-1}. A similar argument proves
equation \eqref{II-2}, and this finishes the proof.
\end{proof}

\begin{rems}
\begin{enumerate}[(1)]
\item We can study the case of special corings
such as a corings constructed by an entwined structures (see
\cite{Brzezinski/Wisbauer:2003}) and other kind of corings, all of them within the bicategory
$\Sf{REM(Bim)}$. This will lead in particular to the study of the
cases of graded modules, comodules over coalgebras over fields,
Hopf modules... etc.

\item The case of the $2$-category whose $0$-cells are all
Grothendieck categories, and whose Hom-Categories consists of
categories of continuous functors i.e. right exact functors which
preserve direct sums (see \cite{Kaoutit:2006} for elementary
details) could be of special interest since it gives another point
of view for the treatments of the problem of equivalence theorems
stated in \cite{Castano/Gomez:1996}, \cite{Castano/Gomez:1998},
\cite{Berbec:2003}, \cite{Chifan/Dascalescu/Nastasescu:2005}.

\item The non unital case, that is, the case concerning the
bicategory of unital bimodules over rings with local units (see \cite{Kaoutit:2006}), leads to the study of
Morita contexts between rings with local units as was developed
in \cite{Abrams:1983} and \cite{Anh/Marki:1987}.
\end{enumerate}
\end{rems}

\providecommand{\bysame}{\leavevmode\hbox
to3em{\hrulefill}\thinspace}
\providecommand{\MR}{\relax\ifhmode\unskip\space\fi MR }
\providecommand{\MRhref}[2]{
} \providecommand{\href}[2]{#2}

\end{document}